\documentclass[10pt,draft,twoside]{amsart}
\usepackage{amssymb}
\usepackage{latexsym}
\usepackage{amsfonts}
\usepackage{amsmath}

\newcommand{\K}{\mathcal K}
\newcommand{\T}{\mathcal T}

\newcommand{\E}{\mathcal E}
\newcommand{\M}{\mathcal M}
\newcommand{\X}{\mathcal X}
\newcommand{\Y}{\mathcal Y}
\newcommand{\cd}[2]{{\sf CD}_{\rm #1}(#2)}

\newcommand{\parti}[2]{{\mathbb P}_{#1}\left(#2\right)}

\newtheorem{theorem}{Theorem}[section]
\newtheorem{proposition}[theorem]{Proposition}
\newtheorem{lemma}[theorem]{Lemma}

\theoremstyle{definition}

\newtheorem{definition}[theorem]{Definition}
\newtheorem{example}[theorem]{Example}

\newcommand{\alt}[1]{{\sf A}_{#1}}

\newcommand{\mat}[1]{{\sf M}_{#1}}
\newcommand{\sy}[1]{{\sf S}_{#1}}

\renewcommand{\sp}[2]{{\sf Sp}_{#1}(#2)}
\newcommand{\pomegap}[2]{{\sf P}\Omega^+_{#1}(#2)}

\newcommand{\sym}[1]{{\sf Sym}\,#1}
\newcommand{\supp}{\operatorname{Supp}}
\newcommand{\diag}{\operatorname{\sf Diag}}

\renewcommand{\wr}{\,\textsf{wr}\,}

\newcommand{\dih}[1]{{\sf D}_{#1}}

\newcommand{\aut}[1]{{\sf Aut}\,{#1}}

\newcommand{\cent}[2]{{\mathbb C}_{#1}(#2)}
\newcommand{\norm}[2]{{\mathbb N}_{#1}\left(#2\right)}

\newcommand{\psl}[2]{\mbox{\sf PSL}_{#1}(#2)}

\renewcommand{\sp}[2]{\mbox{\sf Sp}_{#1}(#2)}

\newcommand{\Z}{\mathbb Z}

\renewcommand{\leq}{\leqslant}
\renewcommand{\geq}{\geqslant}

\begin{document}

\title[Inclusions of innately
transitive groups into wreath products]{Three types of inclusions of innately
transitive\\ permutation groups into wreath products in product action}
\author{Cheryl E. Praeger and Csaba Schneider}
\address[Praeger]{School of Mathematics and Statistics\\
The University of Western Australia\\
35 Stirling Highway, Crawley\\
Western Australia 6009}
\address[Schneider]{Informatics Laboratory\\ Computer and Automation Research Institute of the Hungarian Academy of Sciences\\ 1518 Budapest, Pf.\ 63.}

\email{praeger@maths.uwa.edu.au,
csaba.schneider@sztaki.hu\protect{\newline} {\it WWW:}
www.maths.uwa.edu.au/$\sim$praeger, www.sztaki.hu/$\sim$schneider}

\begin{abstract}
A permutation group is innately transitive if it has a transitive
minimal normal subgroup, and this subgroup is called a plinth. 
In this paper we study three special types of inclusions of innately transitive
permutation groups in wreath products in product action. 
This is achieved by studying the natural Cartesian
decomposition of the underlying set that correspond to the product
action of a wreath product. Previously we identified six classes of
Cartesian decompositions that can be acted upon transitively by an
innately transitive group with a non-abelian plinth. The inclusions
studied in this paper correspond to three of the six classes. We find
that in each case the isomorphism type of the acting group is
restricted, and some interesting combinatorial structures are
left invariant. We also show how to construct examples of inclusions 
for each type.
\end{abstract}

\thanks{{\it Date:} draft typeset \today
The authors acknowledge the support of the Australian Research Council Discovery Grant DP0209706.}

\maketitle

\section{Introduction}\label{s1}

The results of this paper play a key r\^ole in our program to describe
innately transitive subgroups of wreath products in product
action. A permutation group is said to be {\em innately transitive} if it
has a transitive minimal normal subgroup, called a {\em
plinth}; see~\cite{bp}.

Suppose that $G$ is a finite innately transitive subgroup of
$\sym\Omega$  with plinth $M$. Our
aim is to decide whether $G$ can be contained in a subgroup of
$\sym\Omega$ that is permutationally isomorphic to a wreath product
$W=\sym\Gamma\wr \sy\ell$ 
in such a way that $G$ projects onto a
transitive subgroup of $\sy\ell$. Here the group 
$W$ is considered as a permutation group acting on
$\Gamma^\ell$ in product action. Such problems arise in algebraic
combinatorics where often we are given
a combinatorial structure with a subgroup of its automorphism group;
our task is to determine a
larger subgroup of the automorphism group, or, where possible, the full
automorphism group itself. The case where the given group preserves  
additional structure on points, such as a Cartesian decomposition (as 
studied in this paper), is often 
difficult to identify as its existence may not be apparent from the given 
combinatorial information.

If $G$ is contained in a wreath product $W$ as above then the underlying
set $\Omega$ can be identified with the Cartesian product
$\Gamma^\ell$, such that the groups $G$ and $W$ preserve the natural
Cartesian decomposition of $\Gamma^\ell$ (see
Section~\ref{sec2}). Moreover, the permutation representation of $G$ induced
by the natural projection $W\rightarrow \sy\ell$ is equivalent to the
$G$-action on this Cartesian decomposition. 
In~\cite[6-Class~Theorem]{transcs} we identified six classes of Cartesian decompositions 
acted upon transitively by an innately transitive group with a
non-abelian plinth. The names of these classes are $\cd {S}G$, $\cd
1{G}$, $\cd {1S}G$, $\cd {2\sim}G$, $\cd {2\not\sim}G$, and $\cd
3G$ (see Section~\ref{sec2}).
A $G$-invariant Cartesian decomposition of $\Omega$ in a particular
class  leads to 
a special type of  embedding of $G$ into a wreath product in product
action. Cartesian decompositions in $\cd SG$ and $\cd 1G$ were described
in~\cite{transcs}, while $\cd 3G$ was studied
in~\cite{design}. The aim of this paper is to investigate the remaining
three classes, namely $\cd {1S}G$, $\cd{2\sim}G$, and
$\cd{2\not\sim}G$. 

We believe that the classes $\cd {1S}G$, $\cd{2\sim}G$, 
$\cd{2\not\sim}G$, and $\cd 3G$ are the most challenging ones of the 6-Class
Theorem. While the classes $\cd SG$, $\cd 1G$ can be viewed as
natural, the remaining classes correspond to exceptional embeddings of innately
transitive groups into wreath products. The aim of the research
presented here is to understand the exceptional embeddings
that correspond to a Cartesian decomposition in  $\cd {1S}G$, $\cd{2\sim}G$, 
or $\cd{2\not\sim}G$. 
We describe these in as much detail as feasible in our framework. 
Thus this paper contains three main
results that are proved in 
Sections~\ref{7}, \ref{6.5}, and~\ref{6}. As we do not want to litter the
introduction with complicated notation, instead of explicitly 
stating our main results here, we present the following  schema on which 
Theorems~\ref{main2c}, \ref{main2cnot}, and~\ref{main1s} are built.

Let $G$ be an innately transitive group with plinth $M$, where $M$ is 
isomorphic to the direct power of a non-abelian, finite simple group $T$. Suppose that 
$\E\in\cd{2\sim}G\cup\cd{2\not\sim}G\cup\cd{1S}G$. Then in each case we prove
three~properties of the permutation group $G$ and the Cartesian 
decomposition $\E$.

{\bf [1] (Quotient Action Property)}
We study $G$ via its action on a $G$-invariant partition $\overline\Omega$ 
of $\Omega$. We construct a Cartesian decomposition $\overline\E$ of $\overline\Omega$
which is invariant under the action of $\overline G$ and has  characteristics similar to 
those of $\E$. If $\E\in\cd{2\sim}G\cup\cd{1S}G$ then a block in this partition
will have size at most $2^k$, where $k$ depends on $M$, and often $k=0$. 
This last 
statement will not, in general, be true for $\cd{2\not\sim}G$.

{\bf [2] (Factorisation Property)}
We prove that $T$
will admit some special type of factorisation. If
$\E\in\cd{2\sim}G\cup\cd{1S}G$ then this will enable us to severely
restrict the isomorphism type of $T$  (see the {\bf Isomorphism Property} below). 
If $\E\in\cd{2\not\sim}G$ then we will also exclude some isomorphism types 
for $T$.

{\bf [3] (Combinatorial Property)}
We prove that a point stabiliser $G_\omega$ preserves some combinatorial
structure determined by $G$ and $\E$, such as a graph or a generalised graph.

If $\E\in\cd{2\sim}G\cup\cd{1S}G$ then we will also prove an isomorphism property.

{\bf [4] (Isomorphism Property)}
If $\E\in\cd{2\sim}G\cup\cd{1S}G$ then the Factorisation Property is so 
strong that, up to an elementary abelian 2-group,  
we can pinpoint the permutational isomorphism type
of the group $G$. 

Theorems~\ref{main2c}, \ref{main2cnot}, and~\ref{main1s} will be 
built on the above schema.
The converse of these theorems will also hold in the following sense. 
If an innately transitive  permutation group is given together 
with a factorisation
and a combinatorial structure prescribed by the Factorisation Property and
the Combinatorial Property, then we will show how to construct a Cartesian
decomposition $\E$ that belongs to the corresponding class; see
Sections~\ref{7.2}, \ref{7.1}, and~\ref{c1s}. For technical reasons in these
constructions we will require that a point stabiliser in the plinth satisfies
some extra condition to ensure that the partition in the 
Quotient Action Property will contain only trivial blocks. In particular, the 
constructions demonstrate that each instance where the Factorisation Property,
and the Combinatorial Property is satisfied can be realised by an innately 
transitive group $G$ with Cartesian decomposition $\E$ of the appropriate type.

In Section~\ref{sec2} we collect necessary background information on Cartesian 
decompositions and Cartesian systems following the treatment in~\cite{recog}, 
and~\cite{transcs}. Section~\ref{tools} contains some easy lemmas that we need 
for our main theorems. For an innately transitive group $G$, the
Cartesian decompositions in $\E\in\cd{2\sim}G\cup\cd{2\not\sim}G\cup\cd{1S}G$
are studied via their quotient actions, and the required material is 
presented in Section~\ref{quotsec}. Sections~\ref{7}, \ref{6.5}, and~\ref{6} 
are devoted to the classes $\cd{2\sim}G$, $\cd{2\not\sim}G$, and $\cd{1S}G$, 
respectively. We state and prove our main theorems in these three sections.

In this
paper we use the following notation. Permutations act on the right: if $\pi$ is
a permutation and
$\omega$ is a point then the image of $\omega$ under $\pi$ is denoted
$\omega\pi$.   If $G$ is a group acting on a
set $\Omega$ then $G^\Omega$ denotes the subgroup of $\sym\Omega$
induced by $G$.

\section{Cartesian decompositions and Cartesian systems}\label{sec2}

A {\em Cartesian decomposition} of a set $\Omega$ is a set
$\{\Gamma_1,\ldots,\Gamma_\ell\}$ of proper partitions of $\Omega$ such that 
$$
|\gamma_1\cap\cdots\cap\gamma_\ell|=1\quad\mbox{for
all}\quad\gamma_1\in\Gamma_1,\ldots,\gamma_\ell\in\Gamma_\ell.
$$
This property implies that the map
$\omega\mapsto(\gamma_1,\ldots,\gamma_\ell)$, where for
$i=1,\ldots,\ell$ the block $\gamma_i\in\Gamma_i$ is chosen so that
$\omega\in\gamma_i$, is a well defined bijection between $\Omega$ and
$\Gamma_1\times\cdots\times \Gamma_\ell$. Thus the set $\Omega$ can
naturally be
identified with the Cartesian product
$\Gamma_1\times\cdots\times\Gamma_\ell$. 

If $G$ is a permutation group acting on $\Omega$, then a Cartesian
decomposition $\E$ of $\Omega$ is said to be $G$-invariant, if the partitions in
$\E$ are permuted by $G$.
For a permutation group $G$ acting on $\Omega$, the symbol $\cd {}G$
denotes the set of $G$-invariant Cartesian decompositions of
$\Omega$. If $\E\in\cd {}G$ and $G$ acts on $\E$
transitively, then $\E$ is said to be a {\em transitive} $G$-invariant
Cartesian decomposition. The set of transitive $G$-invariant
Cartesian decompositions of $\Omega$ is denoted by $\cd{\rm
  tr}G$. The concept of a Cartesian decomposition was introduced by L.\ G.\ 
Kov\'acs in \cite{kov:decomp} where it is called a system of product imprimitivity.
Kov\'acs suggested that studying $\cd {tr}G$ (using our terminology), for finite 
primitive permutation groups $G$, was the appropriate way to identify 
inclusions of $G$ in wreath products in their product action. His papers 
\cite{kov:decomp} and~\cite{kov:blowups} inspired our work.

Suppose that $G$ is an innately transitive
subgroup of $\sym\Omega$ with plinth $M$, and that $\E$ is a
$G$-invariant Cartesian decomposition of $\Omega$. 
Then we proved
in~\cite[Proposition~2.1]{recog} that each of the $\Gamma_i$ is an
$M$-invariant partition of $\Omega$. 
Choose an element $\omega$ of $\Omega$ and let 
$\gamma_1\in\Gamma_1,\ldots,\gamma_\ell\in\Gamma_\ell$ be such that
$\{\omega\}=\gamma_1\cap\cdots\cap\gamma_\ell$; set $K_i=M_{\gamma_i}$. 
Then
\cite[Lemmas~2.2 and 2.3]{recog} imply that the set $\K_\omega(\E)=\{K_1,\ldots,K_\ell\}$ is invariant
under conjugation by $G_\omega$, 
\begin{eqnarray}\label{csdef1}
\bigcap_{i=1}^\ell K_i&=&M_\omega\quad\mbox{and}\\
\label{csdef2}K_i\left(\bigcap_{j\neq
i}K_j\right)&=&M\quad \mbox{for all}\quad i\in\{1,\ldots,\ell\}.
\end{eqnarray}

For an arbitrary transitive permutation group $M$ on $\Omega$ and a point $\omega\in\Omega$, a set $\K=\{K_1,\ldots,K_\ell\}$ of proper subgroups
of $M$ is said to be 
a {\em Cartesian system of subgroups with respect to $\omega$} 
for $M$, if~\eqref{csdef1} and \eqref{csdef2}~hold. If
$M$ is an abstract group then a set $\{K_1,\ldots,K_\ell\}$ of
proper subgroups satisfying~\eqref{csdef2} is said to be a {\em Cartesian
system}.

\begin{theorem}[Theorem~1.4 and Lemma~2.3~\cite{recog}]\label{bij}
Let $G\leq \sym\Omega$ be an innately transitive permutation group with
plinth $M$. For a fixed $\omega\in\Omega$ 
the correspondence $\E\mapsto\K_\omega(\E)$ is a bijection between the
set of $G$--invariant Cartesian decompositions of $\Omega$
and the set of $G_\omega$--invariant
Cartesian systems of subgroups of $M$ with respect to
$\omega$. Moreover the $G_\omega$--actions on $\E$ and on $\K_\omega(\E)$
are equivalent.
\end{theorem}

Suppose that $G\leq\sym\Omega$ is an innately transitive group with
plinth $M$, and let $\omega\in\Omega$ be fixed. Let $\K$ be a
$G_\omega$-invariant  Cartesian system of subgroups of $M$ with
respect to $\omega$. Then the
previous theorem implies that $\K=\K_\omega(\E)$ for some 
$G$-invariant Cartesian decomposition $\E$ of $\Omega$. Moreover, $\E$ 
consists of the $M$-invariant partitions 
$\{(\omega^K)^m\ |\ m\in M\}$
where $K$ runs through the elements of $\K$. 
This Cartesian decomposition
is usually denoted $\E(\K)$.

Using this
theory we were able to describe those innately transitive subgroups 
of wreath products that have a simple plinth. This led to
a classification of transitive simple and almost simple subgroups of
wreath products in product action
(see~\cite[Theorem~1.1]{recog}).

Now we recall a couple of concepts introduced in~\cite{transcs} to
describe subgroups of direct products. 
Suppose that $M=T_1\times\cdots\times T_k$ where the $T_i$ are
groups, and $k\geq1$. For
$I\subseteq\{T_1,\ldots,T_k\}$ the symbol $\sigma_I:M\rightarrow\prod_{T_i\in
I}T_i$ denotes the natural projection map. We also write
$\sigma_i$ for $\sigma_{\{T_i\}}$. 
A subgroup $X$ of $M$ is said to be a {\em strip} if, for each
$i=1,\ldots,k$, either $\sigma_i(X)=1$ or $\sigma_i(X)\cong X$. The set
of $T_i$ such that $\sigma_i(X)\neq 1$ is called the {\em support} of
$X$ and is denoted $\supp X$. If $T_m\in\supp X$ then we also say that
$X$ {\em covers} $T_m$. Two strips $X_1$ and $X_2$ are said to be 
{\em disjoint} if $\supp X_1\cap\supp X_2=\emptyset$. 
A strip $X$ is said to be {\em full} if
$\sigma_i(X)=T_i$ for all $T_i\in\supp X$,
and it is called {\em non-trivial} if
$|\supp X|\geq 2$. A subgroup $K$ of $M$ is said to be {\em subdirect with 
respect to the direct decomposition $T_1\times\cdots\times T_k$} if
$\sigma_i(K)=T_i$ for all $i$. If $M$ is a finite non-abelian 
characteristically simple group,
then a subgroup $K$ is said to be {\em subdirect} if it is subdirect with
respect to the finest direct decomposition of $M$ (that is, as a direct 
product of simple groups).

Let $M=T_1\times\cdots\times T_k$ be a finite non-abelian characteristically
simple group, where $T_1,\ldots,T_k$ are the simple normal subgroups of
$M$, each isomorphic to the same simple group $T$. If $K$ is a subgroup of 
$M$ and $X$ is a strip in $M$ such that $K=X\times\sigma_{\{T_1,\ldots,T_k\}
\setminus\supp X}(K)$ then we say that $X$ is {\em involved} in $K$. 
A strip $X$ is said to be involved in a Cartesian system $\K$ for $M$
if $X$ is involved in some element of $\K$. 
Note that in this case~\cite[Lemma~2.2]{bad:quasi} and~\eqref{csdef2} 
imply that $X$ is involved in a unique element of $\K$.

A non-abelian plinth of an innately transitive group 
$G$ has the form $M=T_1\times\cdots\times T_k$ where
the $T_i$ are  finite, non-abelian, simple groups. Let $\E\in\cd{}G$ and let
$\K_\omega(\E)$ be a corresponding Cartesian system
$\{K_1,\ldots,K_\ell\}$ for $M$. Then equation~\eqref{csdef2} implies that, 
for all $i\leq k$ and $j\leq\ell$,
\begin{equation}\label{simpfact}
\sigma_i(K_j)\left(\bigcap_{j'\neq j}\sigma_i(K_{j'})\right)=T_i.
\end{equation}
In particular this means that if $\sigma_i(K_j)$ is a proper subgroup
of $T_i$ then $\sigma_i(K_{j'})\neq\sigma_i(K_j)$ for all
$j'\in\{1,\ldots,\ell\}\setminus\{j\}$. 
It is thus important to understand the following sets of subgroups:
\begin{equation}\label{f}
\mathcal F_i(\E,M,\omega)=\{\sigma_i(K_j)\ |\ j=1,\ldots,\ell,\ \sigma_i(K_j)
\neq T_i\}.
\end{equation}
From our remarks above, $|\mathcal F_i(\E,M,\omega)|$ is the number of 
indices $j$ such that
$\sigma_i(K_j)\neq T_i$. 
The set $\mathcal F_i(\E,M,\omega)$ is independent of $i$ up to isomorphism, 
in the sense that if $i_1,\ i_2\in\{1,\ldots,k\}$ and $g\in G_\omega$ 
are such that
$T_{i_1}^g=T_{i_2}$ then $\mathcal
F_{i_1}(\E,M,\omega)^g=\{L^g\ |\ L\in\mathcal F_{i_1}(\E,M,\omega)\}=
\mathcal F_{i_2}(\E,M,\omega)$. This argument
also shows that the subgroups in $\mathcal F_{i_1}(\E,M,\omega)$ are actually 
$G_\omega$-conjugate to the subgroups in $\mathcal F_{i_2}(\E,M,\omega)$.

The following theorem was proved in~\cite[Theorems~5.1 and~6.1]{transcs}. 

\begin{theorem}\label{stripth}
Suppose that $G$ is an innately transitive permutation group with
a non-abelian plinth $M=T_1\times\cdots\times T_k$, where $k\geq1$ and
$T_1,\ldots,T_k$ are
pairwise isomorphic finite simple groups. Let $\E\in\cd {\rm tr}G$
with a corresponding Cartesian system $\K$ for $M$ with respect to a 
point $\omega\in\Omega$, and, for $i=1,\ldots,k$, let
$\mathcal F_i=\mathcal F(\E,M,\omega)$ be defined as in~\eqref{f}.
Then the following all hold.
\begin{enumerate}
\item[(a)] The number $|\mathcal F_i|$ is independent of $i$ and
$|\mathcal F_i|\leq 3$. 
\item[(b)] Suppose that there is a non-trivial, full strip involved in
$\K$. Then $k\geq 2$ and $|\mathcal F_i|\in\{0,\ 1\}$.
\item[(c)] If $X$ is a non-trivial, full strip involved in $\K$ and $|\mathcal F_i|=1$ then $|\supp X|=2$. 
\item[(d)] Set
$\mathcal P=\{\supp X\ |\ X\mbox{ is a non-trivial, full strip
involved in $\K$}\}$. If $\mathcal P\neq\emptyset$ then $\mathcal P$
is a $G$-invariant partition of $\{T_1,\ldots,T_k\}$. Thus if $X_1$ and $X_2$ 
are non-trivial, full strips involved in $\K$ then they are disjoint.
\end{enumerate}
\end{theorem}

The set $\cd{\rm tr}G$ is further subdivided according to
the structure of the subgroups in the corresponding Cartesian systems
as follows. The sets $\mathcal F_i=\mathcal F_i(\E,M,\omega)$ are
defined as in~\eqref{f}. 
\begin{eqnarray*}
\cd{\rm S}G&=&\{\E\in\cd{\rm tr}G\ |\ \mbox{the elements of $\K_\omega(\E)$
are subdirect subgroups in $M$}\};\\
\cd 1G&=&\{\E\in\cd{\rm tr}G\ |\ |\mathcal F_i|=1\mbox{ and $\K_\omega(\E)$
involves no non-trivial, full strip}\};\\
\cd {\rm 1S}G&=&\{\E\in\cd{\rm tr}G\ |\ |\mathcal F_i|=1\mbox{ and $\K_\omega(\E)$
involves non-trivial, full strips}\};\\
\cd {2\sim}G&=&\{\E\in\cd{\rm tr}G\ |\ |\mathcal F_i|=2\mbox{ and
the $\mathcal F_i$ contain two $G_\omega$-conjugate subgroups}\};\\
\cd{2\not\sim}G&=&\{\E\in\cd{\rm tr}G\ |\ |\mathcal F_i|=2\mbox{ and
the subgroups in $\mathcal F_i$ are not $G_\omega$-conjugate}\};\\
\cd 3G&=&\{\E\in\cd{\rm tr}G\ |\ |\mathcal F_i|=3\}.
\end{eqnarray*}

At first glance, it seems that the definitions of the classes $\cd SG$, 
$\cd 1G$, $\cd{1S}G$, $\cd {2\sim}G$, $\cd {2\not\sim}G$, and $\cd 3G$ 
may depend on the choice of the Cartesian system, and hence on the choice 
of the point $\omega$. The
following result, proved in~\cite[Theorems~6.2 and~6.3]{transcs}, 
shows that this is not the case, and shows also that these classes
form a partition of $\cd{\rm tr}G$.  A permutation group is called 
{\em quasiprimitive} if all of its
non-trivial normal subgroups are transitive. A finite quasiprimitive 
group is said to
have {\em compound diagonal type} if it has a unique minimal normal
subgroup $M$, which is non-abelian, and a point stabiliser $M_\omega$
is a non-simple, subdirect subgroup of
$M$. See~\cite{prae:quasi,bad:quasi} for more details..

\begin{theorem}[6-class Theorem]\label{5class}
If $G$ is a finite, innately transitive permutation group with a non-abelian
plinth $M$, then the classes $\cd 1G$, $\cd SG$, $\cd {1S}G$,
$\cd{2\sim}G$, $\cd {2\not\sim}G$, and $\cd 3G$
are independent of the choice of the point $\omega$ used in their definition.
They form a partition of $\cd{\rm tr}G$, and moreover, if $M$ is simple,
then $\cd{\rm tr}G=\cd {2\sim}G$. 
\begin{enumerate}
\item[(a)] If $\cd SG\neq\emptyset$, then $G$ is a quasiprimitive group
of compound diagonal type.
\item[(b)] If $\cd {1S}G\cup\cd {2\sim}G\neq\emptyset$,
then $T$ and the subgroups of the
$\mathcal F_i$ are given by one of the columns of Table~$\ref{newt}$.
\item[(c)] If $\cd {2\not\sim}G\neq\emptyset$, then $T$ admits a
factorisation $T=AB$ with $A, B$ proper subgroups.
\item[(d)] If $\cd 3G\neq\emptyset$, then $T$ is isomorphic to one of the 
groups $\sp{4a}2$ with $a\geq 2$, $\pomegap 83$, or $\sp 62$, and, for
each $i$, the 
subgroups of $\mathcal F_i$ form a strong multiple factorisation of
$T_i$ (see~\cite[Table~V]{bad:fact}), and hence are known explicitly.
\end{enumerate}
\end{theorem}
\begin{table}
$$
\begin{array}{|c|c|c|c|c|}
\hline
T & \alt 6 & \mat{12} & \pomegap 8q & \sp 4{2^a},\ a\geq 2\\
\hline
\mbox{subgroups in $\mathcal F_i$} & \alt 5 & \mat{11} & \Omega_7(q) &
\sp{2}{2^{2a}}\cdot 2\\
\hline
\end{array}
$$
\caption{Table for Theorem~\ref{5class}}
\label{newt}
\end{table}

\section{Toolbox}\label{tools}

In this section we collect the tools, in addition to those
in~\cite{recog,transcs,charfact}, that are necessary for our
investigation of the Cartesian decompositions in $\cd {1S}G$, $\cd
{2\sim}G$, and $\cd{2\not\sim}G$. First
we recall a couple of concepts in graph theory, and then we prove some
group theoretic lemmas.

We introduce the combinatorial structures that are necessary 
for the investigation of the
elements in $\cd {2\sim}G\cup\cd{2\not\sim}G$. 

\begin{definition}
A {\em generalised di-graph} $\Gamma$ is a 4-tuple $(V,E,\beta,\varepsilon)$,
where, $V$ and $E$ are disjoint sets with $V$ non-empty, and
$\beta,\ \varepsilon:E\rightarrow V$ are maps such that $\beta(v)\neq \varepsilon(v)$ for all $v\in V$. 
The elements of $V$ are the {\em vertices}, and the elements of $E$ are
the {\em arcs} of $\Gamma$. If $e\in E$ then $\beta(e)$
is the {\em initial vertex} of $e$, and $\varepsilon(e)$ is the
{\em terminal vertex} of $e$.
A permutation $\alpha\in \sym V\times \sym
E\leq\sym(V\cup E)$ is an {\em automorphism} of $\Gamma$ if
$\alpha(\beta(e))=\beta(\alpha(e))$ and
$\varepsilon(\alpha(e))=\alpha(\varepsilon(e))$ for all $e\in E$. 
\end{definition}

Next we introduce the undirected version of this concept.

\begin{definition}
A {\em generalised graph} $\Gamma$ is a triple $(V,E,\varepsilon)$ where $V$ 
and $E$ are disjoint sets with $V$ non-empty and
$$
\varepsilon:E\rightarrow V^{\{2\}}=\{\{v_1,v_2\}\ |\ v_1,\ v_2\in V,\
v_1\neq v_2\}$$
is a map.  
The elements of $V$ are the {\em vertices} and the elements of $E$ are
the {\em edges} of $\Gamma$. If $e\in E$ then the two elements of 
$\varepsilon(e)$ are said to be {\em adjacent} to $e$.
A permutation $\alpha\in\sym V\times \sym
E\leq\sym(V\cup E)$ is an {\em automorphism} of $\Gamma$ if
$\varepsilon(\alpha(e))=\alpha(\varepsilon(e))$ for all $e\in E$. 
\end{definition}

For the purposes of this paper, a {\em graph} is a generalised graph 
$(V,E,\varepsilon)$ for which $E\subseteq V^{\{2\}}$ and $\varepsilon$ is 
the inclusion map. We usually write this graph simply as $(V,E)$, and with 
the terminology  above, an edge $e=\{v_1,v_2\}$ is adjacent to
$v_1$ and $v_2$ (and vice versa). We will also say that $v_1$ and $v_2$ are 
{\em connected}. 
A graph $(V,E)$ is said to be {\em bipartite} if $V$ has
two non-empty subsets $V_1,\ V_2$ such that $V_1\cap V_2=\emptyset$, $V_1\cup
V_2=V$, and there is no edge between two elements of $V_1$ or between two 
elements of $V_2$. 
The pair $V_1,\ V_2$ is said to be a {\em bipartition} of the graph.

If $\Gamma=(V,E)$ is a graph then the {\em valency} of a vertex $v$ is
defined as the number $|\{v'\in V\ |\ \{v,v'\}\in E\}|$. A graph is said
to be {\em regular} if all vertices have the same valency. A bipartite
graph with a given bipartition is
said to be {\em semiregular} if all vertices in the same part of the
bipartition have the same valency.

A generalised graph can (and will) be viewed as a graph if there is at 
most one edge between any
two vertices. For $n\geq 1$, the complete graph $K_n$ is defined as the graph
in which there are $n$ vertices and any two vertices are connected.

Now we prove some lemmas that are necessary for our investigation.

\begin{lemma}\label{grlem}
Let $\Gamma=(V,E,\varepsilon)$ be a generalised graph such that $E$ is 
non-empty and $\aut\Gamma$ induces a $2$-transitive group on $E$. Then 
either $|V|=2$ or $\Gamma$ is a graph.  In addition, if $\Gamma$ is a graph and 
$\aut\Gamma$ induces a transitive group on $V$, then $\Gamma$ is isomorphic 
to the complete graph $K_3$, or $\Gamma$ is isomorphic to
a vertex-disjoint union of 
$k$ copies of the complete graph $K_2$, for some $k\geq 1$.
\end{lemma}
\begin{proof} Since $E$ is non-empty, we have $|V|\geq2$.
Note that $\aut \Gamma$ must induce a primitive group on $E$. Suppose that
$v_1$ and $v_2$ are vertices of $\Gamma$ such that $v_1$ and $v_2$ are 
connected by some edge in $E$. Then the edges in $E$ that are adjacent to 
$v_1$ and $v_2$ form a block for the action of $\aut \Gamma$ on $E$. Thus 
either $|V|=2$, or 
$v_1$ and $v_2$ are connected by a unique edge in $E$, and so $\Gamma$ is a 
graph.

Assume now that $\Gamma$ is a graph, and, in particular, that $\varepsilon$
is an inclusion map, and that $\aut \Gamma$ is transitive on $V$. This
implies that all 
vertices have the same valency. If this valency is 1 then $\Gamma\cong
k K_2$ for some $k$. So assume that the valency is at least 2 and let 
$v\in V$. Then there
are  edges $e_1$ and $e_2$ such  that  $e_1=\{v,v_1\}$ and  $e_2=\{v,v_2\}$ 
with  $v_1\neq v_2$; in particular $|V|\geq3$.  As  $v_2$ has  valency at
least~2, $v_2$ is  adjacent to an edge  $e_3=\{v_2,v_3\}$ with
$v_3\neq v$. Since $\aut\Gamma$ is 2-transitive on $E$, there is an 
automorphism $\alpha\in\aut\Gamma$ such that $e_1^\alpha=e_1$   and
$e_2^\alpha=e_3$. Thus 
\[
\{v\}^\alpha=(e_1\cap e_2)^\alpha = e_1^\alpha
\cap e_2^\alpha=e_1\cap e_3=\{v,v_1\}\cap\{v_2,v_3\}.\]
Since $v\not\in\{v_2,v_3\}$ and $v_1\ne v_2$, it follows that $v^\alpha 
= v_1=v_3$.  Thus the subgraph spanned by $v,\ v_1,\ v_2$ is a
connected  component of  $\Gamma$ and is a complete graph $K_3$. 
If $|V|\geq 4$, there is a vertex $v_4\not\in\{v,v_1,v_2\}$, and as
$\aut\Gamma$ is transitive on $V$, the connected component of $\Gamma$
containing $v_4$ is also isomorphic to $K_3$. Let $e$ be an edge 
adjacent to $v_4$.  Since $\aut\Gamma$ is 2-transitive on $E$, there   is   an   automorphism   $\beta$   
such   that $(e_1,e_2)^\beta=(e_1,e)$.   Arguing as before, $\{v^\beta\}=
e_1\cap e=\{v,v_1\}\cap e$, but this is the empty set, and we have 
a contradiction. Thus $|V|=3$ and $\Gamma\cong K_3$.
\end{proof}

The next result, which will often be used in complicated arguments, 
 is so easy that its proof is omitted.

\begin{lemma}\label{normlemma1}
Let $A$ and $B$ be subgroups of a group $G$, such that $A\lhd B$
and $\norm GA/A$ is abelian. Then $\norm GA\leq\norm GB$.
\end{lemma}

The following result computes the normaliser of a strip in a direct
product. 

\begin{lemma}\label{normstriplemma}
Let $G_1,\ldots, G_k$ be isomorphic groups, $\varphi_i:G_1\rightarrow G_i$
an isomorphism for $i=2,\ldots,k$, $H_1$ a subgroup of $G_1$, and  $H=\{(h,\varphi_2(h),\ldots,\varphi_k(h))\ |\ h\in H_1\}$
a non-trivial
strip in $G_1\times\cdots\times G_k$. Then 
\begin{equation}\label{normstrip}
\norm {G_1\times\cdots\times G_k}H=\left\{(t,c_2\varphi_2(t),\ldots,c_k\varphi_k(t))\ |\ t\in\norm {G_1}{H_1},\ c_i\in\cent
{G_i}{\varphi_i(H_1)}\right\}.
\end{equation}
\end{lemma}
\begin{proof}
Denote the right hand side of equation~\eqref{normstrip} by $N$, set
$G=G_1\times\cdots\times G_k$, and consider an element  
 $(t,c_2\varphi_2(t), \ldots, c_k\varphi_k(t)) \in N$. Then,  for all
$h\in H_1$,
$$
(h,\varphi_2(h),\ldots,\varphi_k(h))^{(t,c_2\varphi_2(t),\ldots,c_k\varphi_k(t))}=(h^{t},\varphi_2(h^t),\ldots,\varphi_k(h^t))\in
H. 
$$
Hence $N\subseteq \norm {G}H$. Let us prove that
the other inclusion also holds. 
Suppose that
$(t_1,\ldots,t_k)\in \norm{G}H$. Then for all $h\in H_1$
we have that 
$$
(h,\varphi_2(h),\ldots,\varphi_k(h))^{(t_1,\ldots,t_k)}=(h^{t_1},\varphi_2(h)^{t_2},\ldots,\varphi_k(h)^{t_k})\in H,
$$
and so $t_1\in \norm {G_1}{H_1}$ and 
$\varphi_i(h^{t_1})=\varphi_i(h)^{t_i}$ for each $i=2,\dots,k$. This amounts to
saying, for each $i\geq 2$, that $h^{t_1}=h^{\varphi_i^{-1}(t_i)}$ for all $h\in 
H_1$, and hence $t_1\varphi_i^{-1}(t_i)^{-1}\in\cent {G_1}{H_1}$. Therefore
$t_i=c_i\varphi_i(t_1)$ for some $c_i\in\cent{G_i}{\varphi_i(H_1)}$ for
all $i=2,\ldots,k$, and
so 
$$
(t_1,\ldots,t_k)=(t_1,c_2\varphi_2(t_1),\ldots,c_k\varphi_k(t_k)).
$$
Thus $\norm {G}H\subseteq N$, as required. 
\end{proof}

Finally, we need one more result concerning factorisations of
finite simple groups.

\begin{lemma}[Lemma~4.2~\cite{intrans}]\label{simpnorm}
Let $T$ be a finite simple group isomorphic to $\sp 4{2^a}$, where
$a\geq 2$, and let $A,\ B$ be proper isomorphic subgroups of $T$ such
that $T=AB$. Then 
$$
\norm T{A\cap B}=\norm T{A'\cap B'}=A\cap B\quad\mbox{and}\quad
\cent T{A\cap B}=\cent T{A'\cap B'}=1.
$$ 
\end{lemma}

\section{Quotient actions of innately transitive groups}
\label{quotsec}

It is well-known that if $H$ is a transitive permutation group on
$\Omega$ then, for a fixed $\omega\in\Omega$, there is a one-to-one
correspondence between the set $\{H_0\ |\ H_\omega\leq
H_0\leq H\}$ of subgroups and the set of $H$-invariant partitions of 
$\Omega$. The
partition  assigned to $H_0$ by this correspondence is denoted
$\parti H{H_0}$, and is given by
\begin{equation}\label{partdef}
\parti H{H_0}=\left\{\left(\omega^{H_0}\right)^h\ |\ h\in H\right\}.
\end{equation}
In particular, the part of $\parti H{H_0}$ containing $\omega$ is
the $H_0$-orbit $\omega^{H_0}$, and $H_0$ is its setwise stabiliser in $H$.
Note that the next result does not assume that $\Omega$ is finite. 

\begin{lemma}\label{quotlem}
Let $G$ be a permutation group on a set $\Omega$ and $M$ a transitive
normal subgroup of $G$. Suppose that for some $\omega\in\Omega$,
$M_\omega\leq M_0\leq M$ 
and $M_0$ is normalised by $G_\omega$. Then
$\parti M{M_0}$ is $G$-invariant, and if $P\in\parti M{M_0}$ such that
$\omega\in P$ then $G_P=M_0G_\omega$. Moreover, $\parti M{M_0}=\parti
G{M_0G_\omega}$. 
\end{lemma}
\begin{proof}
It is clear from its definition that $\parti M{M_0}$ is $M$-invariant. 
As $M$ is transitive, we have $G=MG_\omega$, and so in order to show that 
$\parti M{M_0}$ is $G$-invariant, it suffices to show that $\parti M{M_0}$
is $G_\omega$-invariant. If $g\in G_\omega$ and $m\in M$ then
$$
\left(\omega^{M_0m}\right)^{g}=\omega^{M_0mg}=\omega^{M_0gm^{g}}=\omega^{gM_0m^g}=\omega^{M_0m^g}\in\parti
M{M_0}.
$$
Hence $\parti M{M_0}$ is $G$-invariant. Thus, by the remarks preceding the 
lemma, $\parti M{M_0}=\parti G{X}$ for a unique subgroup $X$ satisfying
$G_\omega\leq X\leq G$, the part $P=\omega^{M_0}$ containing $\omega$ is 
the $X$-orbit $\omega^{X}$, and $X$ is its setwise stabiliser in $G$. Since 
$\parti M{M_0}$ is $G$-invariant, it follows that 
$G_\omega$ fixes $P$ setwise, as does $M_0$, and by assumption $M_0G_\omega=
G_\omega M_0$ is a subgroup of $G$ containing $G_\omega$.
Since $\omega^{G_\omega M_0}=\omega^{M_0}=\omega^X$, it follows from the 
uniqueness of $X$ that $X=G_\omega M_0$.
\end{proof}

Lemma~\ref{quotlem} can be used 
to construct quotient
actions of innately transitive groups. Suppose that $M$ is
a non-abelian, transitive, minimal normal subgroup of a finite permutation
group $G$, acting on $\Omega$. Let $\omega\in\Omega$ and let $\M$ be a
$G$-invariant partition of the minimal normal subgroups of $M$. If, for 
$I\in \M$, $\sigma_I$ denotes the projection of $M$ to the direct
product of the minimal normal subgroups lying in $I$, then $M_\omega\leq
\prod_{I\in \M}\sigma_I(M_\omega)\leq M$, and we define
$$
\parti {}{\M}=\parti M{\prod_{I\in\M}\sigma_I(M_\omega)}.
$$
As  $\sigma_I(M_\omega)^g=\sigma_{I^g}(M_\omega)$ for all
$I\in\mathcal M$ and $g\in G_\omega$, the subgroup 
$\prod_{I}\sigma_I(M_\omega)$ is normalised by $G_\omega$.
Therefore $\parti{}\M$ is an $M$-invariant partition of $\Omega$.

\section{Cartesian decompositions in $\cd{2\sim}G$}\label{7}

In this section we assume that $G$ is an innately transitive group
acting on $\Omega$ with
a non-abelian plinth $M=T_1\times\cdots\times T_k$ where each  of the $T_i$ is
isomorphic to a finite simple group $T$. Set
$\T=\{T_1,\ldots,T_k\}$,
and fix an $\omega\in\Omega$. Let 
$\overline M_\omega=\norm{T_1}{\sigma_1(M_\omega)}\times\cdots\times
\norm{T_k}{\sigma_k(M_\omega)}$ and let $\overline\Omega$ denote the
$M$-invariant partition $\parti{M}{\overline M_\omega}$ of $\Omega$. 
Using Lemma~\ref{quotlem}, it is easy to see that 
$\overline\Omega$ is $G$-invariant.
Let $\overline\omega$ be the block in
$\overline\Omega$ that contains $\omega$. Then $\overline M_\omega 
= M_{\overline\omega}$ and Lemma~\ref{quotlem}
also implies that $G_{\overline\omega}=\overline M_\omega G_\omega$. 
Let $\overline G$ denote the group induced by 
$G$ on $\overline\Omega$, and let $\overline G_{\overline\omega}$ denote the image in $\overline G$ 
of $G_{\overline\omega}$. 

Suppose that $\E\in\cd{2\sim}G$, and for each $i=1,\dots,k$, let $\mathcal F_i
(\E,M,\omega)=\{A_i, B_i\}$.  Let $\Gamma(G,\E)$
be the generalised graph $(\K_\omega(\E),\T,\varepsilon)$ such that, for 
$i=1,\ldots,k$,  $\varepsilon(T_i)=\{K_{j_1},\ K_{j_2}\}$ where
$\sigma_i(K_{j_1})=A_i$ and $\sigma_i(K_{j_2})=B_i$. For
$i=1,\ldots,\ell$, set $\overline
K_i=\sigma_1(K_i)\times\cdots\times\sigma_k(K_i)$, and let
$\overline\K_\omega(\E)=\{\overline K_1,\ldots,\overline
K_\ell\}$.  

The main result of this section is the following theorem.

\begin{theorem}\label{main2c}
Let the groups $G$ and $M$ be as in the first paragraph of this section. 
If $\E\in\cd {2\sim}G$,
then the properties {\sf Prop2\!$\sim$[a]}--{\sf[d]}
below all hold. 
\end{theorem}

\noindent\textbf{\textsf{Prop2\!$\sim$[a]} (Quotient Action Property).}
The group $M$ is faithful on $\overline\Omega$, and so, if $K$ is a subgroup of $M$, then we
identify $K$ with its image under the action on $\overline\Omega$. 
The set $\overline\K_\omega(\E)$ is a $\overline G_{\overline\omega}$-invariant
Cartesian system of subgroups for $M$ with respect to
$\overline\omega$. Moreover, $\E(\overline\K_\omega(\E))\in\cd {2\sim}{\overline G}$. 

\bigskip

\noindent\textbf{\textsf{Prop2\!$\sim$[b]} (Factorisation Property).}
If $i\in\{1,\ldots,k\}$ then
\begin{enumerate}
\item[(i)] $A_i,\ B_i$
are isomorphic proper subgroups of $T_i$;
\item[(ii)] $A_i$ and $B_i$ are
conjugate under $G_\omega$;
\item[(iii)] $A_iB_i=T_i$,  $A_i\cap B_i=\norm{T_i}{\sigma_i(M_{\omega})}$;
\item[(iv)] $\norm {G_\omega}{T_i}=\{g\in G_\omega\ |\ \{A_i,\
B_i\}^g=\{A_i,\ B_i\}\}$.
\end{enumerate}

\bigskip

\noindent\textbf{\textsf{Prop2\!$\sim$[c]} (Combinatorial Property).}
The group $G_\omega$ induces a group of automorphisms of
the generalised graph $\Gamma(G,\E)$, 
which is transitive on both the vertex-set $\K_\omega(\E)$ and the edge-set 
$\T$, where the $G_\omega$-actions are by conjugation. 
Moreover, if,
for some $i\in\{1,\ldots,k\}$,
$\varepsilon(T_i)=\{K_{j_1},K_{j_2}\}$ and $g\in
\norm{G_\omega}{T_i}$,  then $(A_i,B_i)^g=(A_i,B_i)$ if and only if
$(K_{j_1},K_{j_2})^g=(K_{j_1},K_{j_2})$. 

\bigskip

\noindent\textbf{\textsf{Prop2\!$\sim$[d]} (Isomorphism Property).}
The group $T$, the subgroups of $\mathcal F_i(\E,M,\omega)$, 
and $\sigma_i(M_{\overline\omega})$ are as in
Table~\ref{isomfact1}. 
The group $\overline G$ is permutationally isomorphic to a subgroup of
$\aut M$ acting on $\overline\Omega$. In particular, $M$ is the unique minimal
normal subgroup of $\overline G$, and  $\overline G$ is
quasiprimitive. Moreover, if $T$ is as in rows
1--3 of Table~\ref{isomfact1} then $\overline M_\omega=M_\omega$,  
and so $G\cong\overline G$, as permutation groups. 
Otherwise a block in $\overline\Omega$ has size dividing
$2^k$, the kernel $N$ of the action of $G$ on $\overline\Omega$ is an 
elementary abelian 2-group of rank at most $k$, and $\overline G\cong G/N$. 

\bigskip

A converse of Theorem~\ref{main2c} is also true, see Section~\ref{7.2}.
The following proposition will form the basis for the proof of 
Theorem~\ref{main2c}.

\begin{proposition}\label{c2}
Suppose that $G$, $M$, $T$, $\omega$, $\E$, $\mathcal F_i(\E,M,\omega)$  
are as in the first and second paragraphs of this section. 
Then the isomorphism type of $T$ and that of the subgroups in
$\mathcal F_i(\E,M,\omega)$ are as in one of the rows of Table~$\ref{isomfact1}$.
If one of the rows~$1$--$3$ of Table~$\ref{isomfact1}$ is valid 
then 
\begin{equation}\label{strong}
K=\sigma_1(K)\times\cdots\times\sigma_k(K)\quad\mbox{for}\quad K\in\K_\omega(\E),
\end{equation}
 while if row~$4$ is valid then
\begin{equation}\label{weak}
\sigma_1(K)'\times\cdots\times\sigma_k(K)'\leq
K\quad\mbox{and}\quad
\frac{K}{\sigma_1(K)'\times\cdots\times\sigma_k(K)'}\leq\Z_2^k\quad\mbox{for}\quad
K\in\K_\omega(\E).
\end{equation}
\end{proposition}
\begin{center}
\begin{table}[ht]
$$
\begin{array}{|l|c|c|c|}
\hline
 & T & \mbox{subgroups in }\mathcal F_i(\E,M,\omega) & \sigma_i(M_{\overline\omega})\\
\hline
1 & \alt 6&\alt 5 & \dih{10}\\
\hline
2 & \mat{12}&\mat{11} & \psl 2{11}
 \\
\hline
3 & \pomegap 8q &\Omega_7(q) & {\sf G}_2(q)\\
\hline
4 & \sp 4q,\ q\geq 4\mbox{ even}
& \sp 2{q^2}\cdot 2 & \dih{q^2+1}\cdot 2 \\
\hline
\end{array}
$$
\caption{Table for Proposition~\ref{c2}}
\label{isomfact1}
\end{table}
\end{center}
\begin{proof}
For $i=1,\ldots,k$, we have $\mathcal F_i(\E,M,\omega)=\{A_i,B_i\}$, as 
above. Since $G$ acts transitively on $\mathcal T$ by conjugation, and since,
by the definition of $\cd {2\sim}G$, $A_i$ and $B_i$ are $G_\omega$-conjugate,
the subgroups $A_1,\ldots,A_k$ and $B_1,\ldots,B_k$ are
pairwise isomorphic. Also, since $T_1=A_1B_1$ is a factorisation of a
finite simple group with two isomorphic subgroups, it follows 
from~\cite[Lemma~5.2]{recog} that $T$ and $\mathcal F_1(\E,M,\omega)$ 
are as in Table~\ref{isomfact1}. 
Suppose that $\sigma_1(K_j)'\times\cdots\times\sigma_k(K_j)'\not\leq
K_j$, for some $j$. 
Then it follows
from~\cite[Lemma~2.3]{charfact} that there are 
$i_1,\ i_2\in\{1,\ldots,k\}$ such that
\begin{equation}\label{wrong}
\sigma_{i_1}(K_j)'\times\sigma_{i_2}(K_j)'\not\leq \sigma_{\{i_1,i_2\}}(K_j).
\end{equation}
Suppose first that $\sigma_{i_1}(K_j)=T_{i_1}$. Then~\cite[Lemma~4.2]{transcs} 
implies that $K_j$ involves a full strip $X$ covering $T_{i_1}$. 
However, by Theorem~\ref{stripth}, $X$ cannot be a non-trivial strip since
$\E\in\cd {2\sim}G$. Thus $X=T_{i_1}$, and so $T_{i_1}\leq K_j$. This, 
however, implies that
$\sigma_{i_1}(K_j)\leq \sigma_{\{i_1,i_2\}}(K_j)$, and in this case we
must also have $\sigma_{i_2}(K_j)\leq\sigma_{\{i_1,i_2\}}(K_j)$. Therefore 
$\sigma_{i_1}(K_j)\times\sigma_{i_2}(K_j)=\sigma_{\{i_1,i_2\}}(K_j)$
contradicting~\eqref{wrong}.  Hence $\sigma_{i_1}(K_j)<T_{i_1}$, and the 
same argument shows that $\sigma_{i_2}(K_j)<T_{i_2}$.

Since $\E\in\cd {2\sim}G$, there exist $j_1,\ j_2\in\{1,\ldots,\ell\}\setminus\{j\}$ such
that $\sigma_{i_1}(K_{j_1})<T_{i_1}$ and $\sigma_{i_2}(K_{j_2})<T_{i_2}$. 
It follows from~\eqref{csdef2} that $K_j(K_{j_1}\cap K_{j_2})=M$
(where possibly $j_1=j_2$) and
so 
$$
\sigma_{\{i_1,i_2\}}(K_j)\left(\sigma_{\{i_1,i_2\}}(K_{j_1})\cap 
\sigma_{\{i_1,i_2\}}(K_{j_2})\right)=T_{i_1}\times T_{i_2}.
$$
Note that 
$$
\sigma_{\{i_1,i_2\}}(K_{j_1})\cap 
\sigma_{\{i_1,i_2\}}(K_{j_2})\leq \sigma_{i_1}(K_{j_1})\times
\sigma_{i_2}(K_{j_2})
$$
and hence
$$
\sigma_{\{i_1,i_2\}}(K_j)\left(\sigma_{i_1}(K_{j_1})\times 
\sigma_{i_2}(K_{j_2})\right)=T_{i_1}\times T_{i_2}.
$$
By an observation made at the beginning of this proof, 
$\sigma_{i_1}(K_j),\ \sigma_{i_2}(K_j),\
\sigma_{i_1}(K_{j_1}),\ \sigma_{i_2}(K_{j_2})$ are pairwise
isomorphic. Therefore the factorisation in the previous displayed equation is a full
factorisation (see~\cite[Definition~1.1]{charfact}). On the other hand
\eqref{wrong}~holds, and this
contradicts~\cite[Theorem~1.2]{charfact}. Hence the first inequality of
\eqref{weak} holds for all $K\in\mathcal K_\omega(\E)$.
If $T$ is not as in row~4 of
Table~\ref{isomfact1} then the elements of the $\mathcal F_i$  are finite simple groups, and the stronger equation~\eqref{strong} also follows.

Finally if $T$ is as in row~4 of Table~\ref{isomfact1} then
$\sigma_i(K_j)/\sigma_i(K_j)'\cong\Z_2$, and hence 
$$
\frac{K_j}{\sigma_1(K_j)'\times\cdots\times\sigma_k(K_j)'}\leq\frac{\sigma_1(K_j)\times\cdots\times\sigma_k(K_j)}{\sigma_1(K_j)'\times\cdots\times\sigma_k(K_j)'}\cong\Z_2^k.
$$
\end{proof}

Now we prove Theorem~\ref{main2c}. 

\subsection{Proof of Theorem~\ref{main2c}}

\textbf{\textsf{Prop2\!$\sim$[a]}}
As $\E\in\cd {2\sim}G$ we have that, for all $i$, there are two indices
$j$ such that $\sigma_i(K_j)<T_i$. Hence each of the $\overline K_i$ is a 
proper subgroup of $M$. This also shows that
$\sigma_i(M_\omega)$ is a proper subgroup of $T_i$, and so is
$\norm{T_i}{\sigma_i(M_\omega)}$, for all $i\in\{1,\ldots,k\}$. 
Thus no $T_i$ is contained in $\overline M_\omega$, and so $M$
is faithful on $\overline \Omega$. We will therefore identify
each subgroup $K$ of $M$ with its image under the action on
$\overline\Omega$. 
Set $\K=\K_\omega(\E)$ and $\overline\K=\overline\K_{\omega}(\E)$. 

Next we prove that $\overline\K$ is a $\overline
G_{\overline\omega}$-invariant Cartesian system for $M$ with respect
to $\overline\omega$. If one of the
rows 1--3 of Table~\ref{isomfact1} is valid, then, by Proposition~\ref{c2}, 
$\overline\K=\K$ and
$\overline\omega=\{\omega\}$ and there is nothing to prove. Thus we
suppose that row~4 of Table~\ref{isomfact1} is valid.
First we prove that~\eqref{csdef1} holds. Let $i\in\{1,\ldots,k\}$ and 
$\mathcal F_i=\{A_i,B_i\}$. Then  it
follows from~\eqref{weak} that $A_i'\cap B_i'\leq
\sigma_i(M_\omega)\leq A_i\cap B_i$.
As $A_iB_i=A_i'B_i=A_iB_i'=T_i$ but $A_i'B_i'\neq T_i$, 
we obtain that $|A_i\cap B_i:A_i'\cap
B_i'|=2$. Lemma~\ref{simpnorm} implies that $\norm{T_i}{A_i'\cap
B_i'}=\norm{T_i}{A_i\cap B_i}=A_i\cap B_i$. Hence
$\norm{T_i}{\sigma_i(M_\omega)}=A_i\cap B_i$, and so
$$
\overline K_1\cap\cdots\cap\overline K_\ell=\prod_{i=1}^k(A_i\cap
B_i)=\prod_{i=1}^k\norm{T_i}{\sigma_i(M_\omega)}=M_{\overline \omega}
$$
and condition (\ref{csdef1}) is proved.  
Since~\eqref{csdef2} holds for $\K$ and
$K_i\leq \overline K_i$ for all $i$, we obtain that \eqref{csdef2}~holds
for $\overline\K$ as well. 

We claim that $\overline\K$ is invariant under conjugation by
$G_{\overline\omega}$. Let $i\in\{1,\ldots,k\}$, $j\in\{1,\ldots,\ell\}$, and
$g\in G_\omega$. We denote by $i^g$ the integer in $\{1,\ldots,k\}$
that satisfies
$T_i^g=T_{i^g}$. Then $\sigma_i(K_j)^g=\sigma_{i^g}(K_j^g)$. Thus
$$
(\overline K_j)^g=\left(\prod_{i=1}^k\sigma_i(K_j)\right)^g=\prod_{i=1}^k\sigma_{i^g}\left(K_j^g\right)=\prod_{i=1}^k\sigma_i\left(K_j^g\right).
$$
Since $K_j^g\in \K$, it follows that $\overline K_j^g\in\overline\K$. 
Hence  $\overline\K$ is $G_\omega$-invariant.
Lemma~\ref{quotlem} shows that $G_{\overline\omega}=\overline
M_\omega G_\omega$. Since $\overline M_\omega=\overline K_1\cap\cdots\cap 
\overline K_\ell$ preserves $\overline \K$, we
obtain that $\overline\K$ is also
$G_{\overline\omega}$-invariant. Thus $\overline\K$ is a
$\overline G_{\overline\omega}$-invariant Cartesian system of
subgroups for  $M$
with respect to $\overline\omega$. 

It follows from the definition of $\overline\K$ that $\mathcal
F_i(\E,M,\omega)=\mathcal
F_i(\E(\overline\K),M,\overline\omega)=\{A_i,B_i\}$ for all $i\in\{1,\ldots,k\}$. 
Let $g\in G_\omega$ such that $A_1^g=B_1$ and let $\overline g$ denote
the image of $g$ in its action on $\overline\Omega$. Then
$\overline g\in\overline G_{\overline\omega}$ and clearly
$A_1^{\overline g}=B_1$. Thus $\E(\overline\K)\in\cd{2\sim}{\overline G}$. 

\textbf{\textsf{Prop2\!$\sim$[b]}} Let $i\in\{1,\ldots,k\}$ and choose $j_1,\ j_2\in\{1,\ldots,\ell\}$ such that 
$A_i=\sigma_i(K_{j_1})$ and 
$B_i=\sigma_i(K_{j_2})$. 
It is clear by the definition of $\cd{2\sim}G$ 
that {\sf Prop2\!$\sim$[b](i)-(ii)} hold
for $A_i$ and $B_i$. Since $K_{j_1}K_{j_2}=M$ we have that 
$\sigma_i(K_{j_1})\sigma_i(K_{j_2})=\sigma_i(M)$ and so
$A_iB_i=T_i$. We showed in the proof of \textsf{Prop2\!$\sim$[a]} that
$A_i\cap B_i=\norm{T_i}{\sigma_i(M_\omega)}$, and so 
{\sf Prop2\!$\sim$[b](iii)} also holds. Finally, let $g\in
G_\omega$ such that $\{A_i,B_i\}^g=\{A_i,B_i\}$. Since $A_i,\ B_i\leq T_i$, it
follows that $T_i\cap T_i^g\neq 1$, and so $T_i^g=T_i$. Conversely, if
$T_i^g=T_i$ with some $g\in G_\omega$, then $\sigma_i(
K)^g=\sigma_i(K^g)$ for all $K\in\K$. Thus the uniqueness of $\{j_1,\
j_2\}$ yields that $g$ fixes $\{K_{j_1},\  K_{j_2}\}$. Therefore $g$
fixes $\{\sigma_i(K_{j_1}),\ \sigma_i(K_{j_2})\}=\{A_i,B_i\}$. Thus all
properties in {\sf Prop2\!$\sim$[b]} hold.

\textbf{\textsf{Prop2\!$\sim$[c]}}
Let $\Gamma$ denote $\Gamma(G,\E)$.
It follows from the definition of $\Gamma$ that
the action of $G_\omega$ by conjugation
is transitive on the vertex set $\K$ and on the
edge set $\T$ of
$\Gamma$. We claim that $G_\omega$ preserves adjacency in
$\Gamma(G,\E)$. Let $i_1\in\{1,\ldots,k\}$ 
with $\varepsilon(T_{i_1})=\{K_{j_1},K_{j_2}\}$ 
and let $g\in G_\omega$. Let
$i_2\in\{1,\ldots,k\}$ such that $T_{i_1}^g=T_{i_2}$. Then
$\sigma_{i_1}(K_{j_1})^g=\sigma_{i_2}(K_{j_1}^g)$ and
$\sigma_{i_1}(K_{j_2})^g=\sigma_{i_2}(K_{j_2}^g)$. Thus
$\varepsilon(T_{i_1}^g)=\varepsilon(T_{i_2})=\{K_{j_1}^g,K_{j_2}^g\}=\varepsilon(T_{i_1})^g$,
as required.

Suppose that $g\in
\norm{G_{\omega}}{T_i}$. If $g$ is such that $A_i^g=A_i$ and $B_i^g=B_i$ then,
as $g\in \norm{G_{\omega}}{T_i}$,
$A_i=A_i^g=\sigma_i(K_{j_1})^g=\sigma_i(K_{j_1}^g)$. As $j_1$ is the unique
integer in $\{1,\ldots,\ell\}$ such that $\sigma_i(K_{j_1})=A_i$, we
obtain that $K_{j_1}^g=K_{j_1}$, and also $K_{j_2}^g=K_{j_2}$. If $g\in\norm{G_{\omega}}{T_i}$
is such that $K_{j_1}^g=K_{j_1}$ and $K_{j_2}^g=K_{j_2}$ then it also
follows that
$A_i^g=\sigma_i(K_{j_1})^g=\sigma_i(K_{j_1}^g)=\sigma_i(K_{j_1})=A_i$, and,
of course, $B_i^g=B_i$. Thus the compatibility condition between $A_i$, $B_i$,
and $\Gamma(G,\E)$ in {\sf Prop2\!$\sim$[c]} also holds.

\textbf{\textsf{Prop2\!$\sim$[d]}}
It follows from Proposition~\ref{c2} that $T$, and the subgroups
$A_i,\ B_i$ of $\mathcal F_i(\E,M,\omega)$ are as in
Table~\ref{isomfact1}. As $M_{\overline\omega}=\overline
K_1\cap\cdots\cap\overline K_\ell$, we obtain that
$\sigma_i(M_{\overline\omega})=A_i\cap B_i$ for all $i$. Since 
$T_i=A_iB_i$ is a factorisation of $T_i$ with two isomorphic
subgroups, we obtain from the~\cite{atlas} in rows~1--2, from~\cite[3.1.1(vi)]{kleidman}
in row~3, and from~\cite[3.2.1(d)]{lps:max} in row~4 of Table~\ref{isomfact1} that the
$\sigma_i(M_{\overline\omega})$-column of Table~\ref{isomfact1} is
correct.
As $M_{\overline\omega}=\overline M_\omega$ is the direct
product of its projections under the $\sigma_i$, and such a projection is
self-normalising in $T_i$ (by Lemma~\ref{simpnorm}), we obtain that
$M_{\overline\omega}$ is a self-normalising subgroup in
$M$. Thus~\cite[Theorem~4.2A]{dm} implies that
$\cent{\sym\overline\Omega}M=1$. Hence $M$ is the unique minimal normal
subgroup of $\overline G$, and so $\overline G$ can be embedded into a
subgroup of $\aut M$. In particular, $\overline G$ is quasiprimitive. 

Suppose that 
$$
\underline
K_i=\sigma_1(K_i)'\times\cdots\times\sigma_k(K_i)'\quad\mbox{for
all}\quad i\in\{1,\ldots,\ell\}
$$
and set $\underline
M_\omega=\underline K_1\cap\cdots\cap\underline K_\ell$. It follows
from Theorem~\ref{c2} that $\underline K_i\leq K_i\leq\overline K_i$
and that  $\underline K_i= K_i=\overline K_i$ if $T$ is as in one of
the rows 1--3 of Table~\ref{isomfact1}; thus $\underline M_\omega\leq
M_\omega\leq M_{\overline \omega}$ also holds. 
If $T$ is as in one of
the rows 1--3 of Table~\ref{isomfact1}, then 
$\underline
M_\omega=M_\omega=M_{\overline \omega}$. Thus $\overline\Omega$ can be identified with $\Omega$, and so
the groups $G$ and $\overline G$ are permutationally isomorphic.

Suppose now that $T$ is as in row~4 of
Table~\ref{isomfact1}. Then, by~\cite[3.2.1(d)]{lps:max},
$\sigma_i(M_{\overline \omega})\cong A_i\cap B_i\cong\dih{q^2+1}\cdot 2$ and
$\sigma_i(\underline M_\omega)\cong A_i'\cap B_i'\cong\dih{q^2+1}$ for all
$i\in\{1,\ldots,k\}$. It follows from Lemma~\ref{simpnorm} that
$$
\norm{T_i}{\sigma_i(\underline M_\omega)}=\norm{T_i}{\sigma_i(M_{\overline
\omega})}=\sigma_i(M_{\overline \omega})\quad\mbox{for all}\quad i\in\{1,\ldots,k\}.
$$ 
Hence we obtain that
$\norm M{M_\omega}\leq \norm M{\underline M_\omega}=M_{\overline \omega}$. On
the other hand, as $\norm M{\underline M_\omega}/\underline M_\omega$ is
abelian, Lemma~\ref{normlemma1} gives $\norm M{\underline M_\omega}\leq \norm M{M_\omega}$. Thus
$\norm M{\underline M_\omega}=\norm M{M_\omega}$. Therefore $\norm
M{M_\omega}/M_\omega$ is an elementary abelian 2-group of rank at most $k$,
and, by~\cite[Theorem~4.2A]{dm}, a block in $\overline\Omega$ also has size
dividing $2^k$.
Therefore $N$ is also an elementary
abelian $2$-group of rank at most $k$. \hfill$\Box$

\medskip

\subsection{A converse of Theorem~\ref{main2c}}\label{7.2}
Theorem~\ref{main2c} can be reversed in the following sense.
Let $G$ be a finite  innately transitive group on $\Omega$ with a non-abelian plinth
$M$ and let $T_1,\ldots,T_k$ be the simple direct factors of $M$. Assume that, 
for $\omega\in\Omega$, the point stabiliser
$M_\omega$ can be decomposed as
$M_\omega=\sigma_1(M_\omega)\times\cdots\times\sigma_k(M_\omega)$. Set
$\T=\{T_1,\ldots,T_k\}$. 
Suppose that $A_1,\ B_1$ are subgroups of
$T_1$ and $\Gamma=(V,\T,\varepsilon)$ is a generalised graph, such
that properties \textsf{Prop2\!$\sim$[b]} and \textsf{Prop2\!$\sim$[c]} hold.
More precisely,
\begin{enumerate}
\item[(i)] $A_1,\ B_1$
are isomorphic proper subgroups of $T_1$;
\item[(ii)] $A_1$ and $B_1$ are
conjugate under $G_\omega$;
\item[(iii)] $A_1B_1=T_1$,  $A_1\cap B_1=\sigma_1(M_{\omega})$;
\item[(iv)] $\norm {G_\omega}{T_1}=\{g\in G_\omega\ |\ \{A_1,\
B_1\}^g=\{A_1,\ B_1\}\}$.
\end{enumerate}
Assume, moreover, 
that $G_\omega$ induces a vertex and edge-transitive group of automorphisms
of $\Gamma$, where the $G_\omega$-action on $\T$ is by
conjugation, and that, if  $\varepsilon(T_1)=\{v_1,v_2\}$ in $\Gamma$, then
the following holds:
\begin{equation}\label{comp}
\mbox{if $g\in\norm{G_\omega}{T_1}$ then
$(A_1,B_1)^g=(A_1,B_1)$ if and only if
$(v_1,v_2)^g=(v_1,v_2)$}.
\end{equation}
We construct, as follows, a $G$-invariant Cartesian decomposition $\E$ in 
$\cd{2\sim}G$, 
such that $\Gamma\cong\Gamma(G,\E)$ and $\mathcal F_1(\E,M,\omega)=\{A_1,
B_1\}$.

For $i=1,\ldots,k$, choose $g_i\in G_\omega$ such that
$T_1^{g_i}=T_i$. For each element $v\in V$ set $K_v=\prod_{i=1}^k
K_{v,i}$ where, for $i=1,\ldots,k$, the subgroup $K_{v,i}$ is defined 
as follows  (noting that $\varepsilon(T_i)=\{v_1^{g_i},v_2^{g_i}\}$). Set
$K_{v_1^{g_i},i}=A_1^{g_i}$, $K_{v_2^{g_i},i}=B_1^{g_i}$, and
$K_{v,i}=T_i$ for all $v\in V\setminus\{v_1^{g_i},v_2^{g_i}\}$.

We claim that each of the $K_{v,i}$ is well-defined, that is, its definition
is independent of the choice of the $g_i$. Suppose that $g_i,\ g_i'\in
G_\omega$ are such that $T_1^{g_i}=T_1^{g_i'}=T_i$ for some $i$. Note that, as
$G_\omega$ induces a group of automorphisms of $\Gamma$,  in this
case
$\{v_1^{g_i},v_2^{g_i}\}=\varepsilon(T_i)=\{v_1^{g_i'},v_2^{g_i'}\}$. 
Thus if $v\not\in \{v_1^{g_i},v_2^{g_i}\}$ then we would define $K_{v,i}$ as
$T_{i}$ using either $g_i$ or $g_i'$.
Suppose next
that $v_1^{g_i}=v_1^{g_i'}$ and $v_2^{g_i}=v_2^{g_i'}$. 
Then $g_ig_i'{}^{-1}\in\norm{G_\omega}{T_1}\cap
(G_\omega)_{v_1}$ and so, by~\eqref{comp}, $g_ig_i'{}^{-1}\in\norm{G_\omega}
{A_1}\cap\norm{G_\omega}{B_1}$. Thus $A_1^{g_ig_i'{}^{-1}}=A_1$ and 
$B_1^{g_ig_i'{}^{-1}}=B_1$; and so $A_1^{g_i}=A_1^{g_i'}$ and $B_1^{g_i}
=B_1^{g_i'}$. Therefore, using either $g_i$ or $g_i'$, we would define 
$K_{v_1^{g_i},i}$ as $A_1^{g_i}$ and $K_{v_2^{g_i},i}$ as $B_1^{g_i}$.
The other possibility is that $v_1^{g_i}=v_2^{g_i'}$ and
$v_1^{g_i'}=v_2^{g_i}$. Then $g_ig_i'{}^{-1}$ is in $\norm{G_\omega}{T_1}$ 
and interchanges $v_1$ and $v_2$, and so by property (iv) above and condition 
(\ref{comp}), $g_ig_i'{}^{-1}$ also swaps $A_1$ and $B_1$. 
For $v=v_1^{g_i}=v_2^{g_i'}$ we would, using $g_i'$, define $K_{v,i}$ as 
$B_1^{g_i'}=(A_1^{g_ig_i'^{-1}})^{g_i'}=A_1^{g_i}$, and similarly, for 
$v=v_2^{g_i}=v_1^{g_i'}$ we would, using $g_i'$, define $K_{v,i}$ as 
$A_1^{g_i'}=(B_1^{g_ig_i'^{-1}})^{g_i'}=B_1^{g_i}$. Thus the definitions 
of all the $K_{v,i}$ are the same whether we use $g_i$ or $g_i'$.

Let $\K=\{K_v\ |\ v\in V\}$. 
We claim that $\K$ is a $G_\omega$-invariant
Cartesian system for $M$ with respect to $\omega$. 
First note that the $K_v$ are direct products of their projections
and, for all $i$,
$$
\bigcap_{v\in V}K_{v,i}=A_1^{g_i}\cap B_1^{g_i}=(A_1\cap
B_1)^{g_i}=\sigma_1(M_\omega)^{g_i}=\sigma_{i}(M_\omega).
$$
Therefore 
$$
\bigcap_{v\in V}K_v=\prod_{i=1}^k\sigma_i(M_\omega)=M_\omega.
$$
Hence~\eqref{csdef1} holds.
The choice of $A_1$ and $B_1$ is such that $T_1=A_1B_1$, and the definition 
of $K_v=\prod_{i}K_{v,i}$ implies that, for each $i$ and $v$, 
$$
K_{v,i}\left(\bigcap_{v'\neq v} K_{v',i}\right)=T_i.
$$
As $K_{v,i}\leq K_v$ for all $i$ and $v$, it follows that 
$K_v\left(\bigcap_{v'\neq v}K_{v'}\right)=M$ for all
$v$. Thus \eqref{csdef2}~holds and $\K$ is a Cartesian system for
$M$ with respect to $\omega$.  
Now we prove that the set $\K$  is invariant under
conjugation by $G_\omega$. 
Let $v\in V$, $i_1,\ i_2\in\{1,\ldots,k\}$ and $g\in G_\omega$ such
that $T_{i_1}^g=T_{i_2}$. 
We claim that $K_{v,i_1}^g=K_{v^g,i_2}$. Suppose first that
$v=v_1^{g_{i_1}}$. As $g$ induces
an automorphism of $\Gamma$, we obtain that
$v_1^{g_{i_1}g}\in\varepsilon(T_{i_1}^g)=\varepsilon(T_{i_2})=\{v_1^{g_{i_2}},v_2^{g_{i_2}}\}$. 
If $v^g=v_1^{g_{i_1}g}=v_1^{g_{i_2}}$, then 
$g_{i_1}gg_{i_2}^{-1}$ stabilises $(v_1,v_2)$, and hence,
by~\eqref{comp}, normalises $A_1$ and $B_1$, so that $K_{v,i_1}^g=A_1^{g_{i_1}g}=A_1^{g_{i_2}}$, and
$K_{v^g,i_2}=A_1^{g_{i_2}}$. Thus $K_{v,i_1}^g=K_{v^g,i_2}$.
Similar arguments show that $K_{v,i_1}^g=K_{v^g,i_2}$ holds in
all other cases.
Therefore 
$$
K_v^g=\left(\prod_{i=1}^k K_{v,i}\right)^g=\prod_{i=1}^k K_{v^g,i}=K_{v^g}.
$$
Hence $\K$ is $G_\omega$-invariant. We also note that the $G_\omega$-actions on $V$ and on $\K$ are equivalent. Thus $\K$ is a
$G_\omega$-transitive Cartesian system of subgroups in $M$ with
respect to $\omega$, and it follows from the definition of the $K_v$
that $\E(\K)\in\cd {2\sim}G$, $\Gamma\cong\Gamma(G,\E(\K))$, and 
$\mathcal F_1(\E(\K),M,\omega)=\{A_1,B_1\}$. 

\medskip

One aim of this section is to describe those innately transitive 
permutation groups $G$ for which $\cd{2\sim}G$ is non-empty. Our results 
show that, if $\cd {2\sim}G\ne \emptyset$, then the following all hold: 
the isomorphism type of such groups is restricted (see \textsf{Prop2\!
$\sim$[d]}), and a point stabiliser 
in the plinth also satisfies some interesting properties, expressed in 
\textsf{Prop2\!$\sim$[b]}. Moreover, such 
groups $G$ act on a generalised graph (see \textsf{Prop2\!$\sim$[c]}) 
whose edge set is intrinsic to the abstract group theoretic structure 
of $G$. This suggests that the conjugation action of $G$ 
on the simple direct factors of the plinth may, in certain cases, 
predetermine the existence of Cartesian decompositions in 
$\cd{2\sim}G$. This problem would be very interesting to address in more 
detail, but it would distract us from the main focus of 
the present work. We only illustrate this phenomenon with the following 
example 

\begin{example}
Suppose that $G$ is a quasiprimitive permutation group on $\Omega$ with 
a unique minimal normal subgroup $M=T_1\times\cdots\times T_k$, 
where $k\geq 4$ and $T_1,\ldots,T_k$ are finite simple groups, isomorphic 
to one of the groups $T$ in Table~\ref{isomfact1}. Assume further that 
the conjugation action of $G$ induces a 2-transitive 
permutation group on the set $\T=\{T_1,\ldots,T_k\}$. Let $\omega\in\Omega$. 
Then $G=MG_\omega$, and so the $G_\omega$-action on $\T$ is also
2-transitive. Let $\Gamma=(V,\T,\varepsilon)$ be a 
generalised graph such that $G_\omega$ induces a vertex-transitive 
group of automorphisms on $\Gamma$ where the $G_\omega$-action on $\T$ is 
by conjugation. Then Lemma~\ref{grlem} implies that 
$\Gamma$ is isomorphic to the union of 
$k$ copies of the complete graph $K_2$. This shows 
that if $\E\in\cd{2\sim}G$, then $|\E|=2k$. 
Further, if $K$ is a subgroup in the Cartesian system $\K_\omega(\E)$, then
$K$ corresponds to a vertex of $\Gamma$ that is adjacent to a unique edge of
$\Gamma$. Thus there is a unique $i\in\{1,\ldots,k\}$ such that 
$\sigma_i(K)\ne T_i$, and, since $\K_\omega(\E)$ involves no strips,
there is a unique $i$ such that $T_i\not\leq K$. This shows that the corresponding 
embedding of $G$ into the full stabiliser in $\sym\Omega$ of $\E$ is as 
in~\cite[Theorem~1.1(b)]{transcs}. 
\end{example}

\section{Cartesian decompositions in $\cd{2\not\sim}G$}\label{6.5}

In this section we assume that $G$ is a finite innately transitive group
acting on $\Omega$ with
a non-abelian plinth $M=T_1\times\cdots\times T_k$ where each  of the $T_i$ is
isomorphic to a finite simple group $T$. Set
$\T=\{T_1,\ldots,T_k\}$,
and fix an $\omega\in\Omega$. 

Suppose that $\E\in\cd{2\not\sim}G$ and let $\K_\omega(\E)=\{K_1,\ldots,K_\ell\}$ be the corresponding Cartesian system of subgroups. For
$i=1,\ldots,\ell$, set $\overline
K_i=\sigma_1(K_i)\times\cdots\times\sigma_k(K_i)$, and let
$\overline\K_\omega(\E)=\{\overline K_1,\ldots,\overline
K_\ell\}$. Let $\overline M_\omega=\overline K_1\cap\cdots\cap\overline 
K_\ell$, and note that $\overline M_\omega$ is the direct product of its 
projections, that is to say,
$\overline M_\omega=\sigma_1(\overline M_\omega)\times\cdots\times \sigma_k(\overline M_\omega)$. Let $\overline\Omega$ denote the $M$-invariant partition
$\parti M{\overline M_\omega}$ of $\Omega$. 
It is routine to check that 
the conjugation action of
$G_\omega$ permutes the subgroups $\overline K_1,\ldots,\overline K_\ell$, and 
so their intersection $\overline M_\omega$ is invariant under $G_\omega$. 
Thus Lemma~\ref{quotlem} shows that $\overline\Omega$ is $G$-invariant.
Let $\overline\omega$ be the block in
$\overline\Omega$ that contains $\omega$. Then $\overline M_\omega=
M_{\overline\omega}$, and  Lemma~\ref{quotlem}
also implies that $G_{\overline\omega}=\overline M_\omega G_\omega$. 
Let $\overline G$ denote the group induced by $G$ on $\overline\Omega$, so
that $\overline G_{\overline\omega}$ is the image in $\overline G$ of 
$G_{\overline\omega}$.

Define a generalised di-graph $\Gamma(G,\E)=(\K_\omega(\E),\T,\beta,
\varepsilon)$ for the given Cartesian decomposition $\E$ as follows. Let 
$A_1$ and $B_1$ be the subgroups of $T_1$ such that $\mathcal F_1(\E,M,\omega)=\{A_1,B_1\}$. Then 
for each $i$ there are unique indices $j_1$ and $j_2$ such that 
$\sigma_i(K_{j_1})$ is $G_\omega$-conjugate to $A_1$ and 
$\sigma_i(K_{j_2})$ is $G_\omega$-conjugate to $B_1$. Set 
$\beta(T_i)=K_{j_1}$ and $\varepsilon(T_i)=K_{j_2}$, and let $A_i$ 
and $B_i$ denote 
$\sigma_i(K_{j_1})$ and $\sigma_i(K_{j_2})$, respectively. Thus the subgroups
$A_1,\ldots,A_k$ are pairwise $G_\omega$-conjugate, and so are the subgroups
$B_1,\ldots,B_k$. On the other hand, if $i,\ j\in\{1,\ldots,k\}$, then 
$A_i$ is not $G_\omega$-conjugate to $B_j$. 

The main result of this section is the following theorem.

\begin{theorem}\label{main2cnot}
Let the groups $G$ and $M$ be as in the first paragraph of this section. 
If $\E\in\cd {2\not\sim}G$,
then the properties {\sf Prop2\!$\not\sim$[a]}--{\sf Prop2\!$\not\sim$[c]}
below all hold. 
\end{theorem}

\noindent\textbf{\textsf{Prop2\!$\not\sim$[a]} (Quotient Action Property).}
The group $M$ is faithful on $\overline\Omega$, and so, if $K$ is a subgroup of $M$, then we
identify $K$ with its image under the action on $\overline\Omega$. 
The set $\overline\K_\omega(\E)$ is a $\overline G_{\overline\omega}$-invariant
Cartesian system of subgroups for $M$ with respect to
$\overline\omega$. Moreover, $\E(\overline\K_\omega(\E))\in\cd {2\not\sim}{\overline G}$. 

\bigskip

\noindent\textbf{\textsf{Prop2\!$\not\sim$[b]} (Factorisation Property).}
If $i\in\{1,\ldots,k\}$ then
\begin{enumerate}
\item[(i)] $A_i,\ B_i$
are proper subgroups of $T_i$;
\item[(ii)] $A_i$ and $B_i$ are
not conjugate under $G_\omega$;
\item[(iii)] $A_iB_i=T_i$ and $A_i\cap B_i=\sigma_i(\overline M_\omega)$;
\item[(iv)] $\norm {G_\omega}{T_i}=\norm{G_\omega}{A_i}=\norm{G_\omega}{B_i}$.
\end{enumerate}

\bigskip

\noindent\textbf{\textsf{Prop2\!$\not\sim$[c]} (Combinatorial Property).}
The group $G_\omega$ induces a group of automorphisms of
the generalised di-graph $\Gamma(G,\E)$, 
which is transitive on both the vertex-set $\K_\omega(\E)$ and the arc-set 
$\T$, where the $G_\omega$-actions are by conjugation.

\bigskip

The observant reader may notice that our conclusions in this section are
considerably 
weaker than those in Section~\ref{7}, as there is no counterpart of {\sf Prop2$\!\sim$[d]}. The reason for this is simple: in the previous section the finite
simple group $T$ admitted a factorisation with two proper, isomorphic 
subgroups, and so the isomorphism type of $T$, and hence that of $G$, could
be restricted. No such factorisation is guaranteed to exist here. On the
other hand, for some $i$, the subgroups $A_i$ and $B_i$ may be isomorphic even
though they are not $G_\omega$-conjugate. It is easy to see, and is left to 
the reader, that claims similar to those in {\sf Prop2$\!\sim$[d]} 
are valid in this case. We formulate the following related problem.

\medskip

\noindent{\bf Problem.} Let $G$, $M$, and $\omega$ 
be as in the first paragraph of this section
and let $\E\in\cd{2\not\sim}G$ such that $\mathcal F_1(\E,M,\omega)$ contains
two isomorphic subgroups. Is it always true that there is an innately 
transitive subgroup $H$ of $\sym\Omega$, having the same plinth $M$ as $G$, 
such that $\E\in\cd{2\sim}H$?

Next we prove Theorem~\ref{main2cnot}.

\subsection{Proof of Theorem~\ref{main2cnot}}
\textbf{\textsf{Prop2\!$\not\sim$[a]}}\quad
As $\E\in\cd {2\not\sim}G$ we have that, for all $i$, there are two indices
$j$ such that $\sigma_i(K_j)<T_i$. Thus each of the $\overline K_i$ is a 
proper subgroup of $M$. This also shows that
$\sigma_i(\overline M_{\omega})$ is a proper subgroup of $T_i$ for all $i\in\{1,\ldots,k\}$. 
Thus no $T_i$ is a subgroup of $\overline M_{\omega}$, and so $M$ must be
faithful on $\overline \Omega$. 
Set $\K=\K_\omega(\E)$ and $\overline\K=\overline\K_\omega(\E)$.

Next we prove that $\overline\K$ is a $\overline
G_{\overline\omega}$-invariant Cartesian system for $M$ with respect
to $\overline\omega$. 
Equation~\eqref{csdef1} holds because of the 
definition of $\overline M_\omega=M_{\overline\omega}$. 
Since~\eqref{csdef2} holds for $\K$ and
$K_i\leq \overline K_i$ for all $i$, we obtain that \eqref{csdef2}~holds
for $\overline\K$ as well. 
We claim that $\overline\K$ is invariant under conjugation by
$G_\omega$. Let $i\in\{1,\ldots,k\}$, $j\in\{1,\ldots,\ell\}$, and
$g\in G_\omega$. We denote by $i^g$ the integer in $\{1,\ldots,k\}$
that satisfies
$T_i^g=T_{i^g}$. Then $\sigma_i(K_j)^g=\sigma_{i^g}(K_j^g)$. Thus
$$
(\overline K_j)^g=\left(\prod_{i=1}^k\sigma_i(K_j)\right)^g=\prod_{i=1}^k\sigma_{i^g}\left(K_j^g\right)=\prod_{i=1}^k\sigma_i\left(K_j^g\right).
$$
Since $K_j^g\in \K$, it follows that $(\overline K_j)^g\in\overline\K$. 
Lemma~\ref{quotlem} shows that $G_{\overline\omega}=\overline
M_\omega G_\omega$. Then, since $\overline M_\omega=\overline K_1\cap
\cdots\cap \overline K_\ell$ preserves $\overline \K$, and $\overline\K$ 
is $G_\omega$-invariant, we obtain that $\overline\K$ is also
$G_{\overline\omega}$-invariant. Thus $\overline\K$ is a
$\overline G_{\overline\omega}$-invariant Cartesian system of subgroups for $M$
with respect to $\overline\omega$. 

It follows from the definition of $\overline\K$ that $\mathcal
F_i(\E,M,\omega)=\mathcal
F_i(\E(\overline\K),M,\overline\omega)=\{A_i,B_i\}$ for all
$i\in\{1,\ldots,k\}$.  Let $\overline g\in \overline
G_{\overline\omega}$ such that $A_1^{\overline g}=B_1$ and let $g$ denote its
preimage in $G_{\overline\omega}$.  Then $g=mg_1$ for some $m\in
\overline M_\omega$ and $g_1\in G_\omega$. As $\overline M_\omega$ is
the intersection of the $\overline K_j$, we obtain that 
$$
\sigma_1(m)\in\sigma_1\left(\overline K_1\cap\cdots\cap\overline K_\ell\right)\leq
\sigma_1\left(\overline K_1\right)\cap\cdots\cap\sigma_1\left(\overline K_\ell\right)=
\sigma_1(K_1)\cap\cdots\cap\sigma_1(K_\ell)=A_1\cap B_1.
$$
Therefore $A_1^m=A_1$, and so $A_1^{g_1}=B_1$. As $g_1\in G_\omega$
and $\E\in \cd{2\not\sim}G$, this is a contradiction, and so $A_1$ and
$B_1$ are not $\overline G_{\overline\omega}$-conjugate.  Thus
$\E(\overline\K)\in\cd{2\not\sim}{\overline G}$.

\textbf{\textsf{Prop2\!$\not\sim$[b]}}\quad Let $i\in\{1,\ldots,k\}$ and 
let $j_1,\ j_2\in\{1,\ldots,\ell\}$ be such that $A_i=\sigma_i(K_{j_1})$ and 
$B_i=\sigma_i(K_{j_2})$. By the definition of $\cd {2\not\sim}G$, it
follows that {\sf Prop2\!$\not\sim$[b]}(i)--(ii) hold
for $A_i$ and $B_i$. Since $K_{j_1}K_{j_2}=M$ we have that 
$\sigma_i(K_{j_1})\sigma_i(K_{j_2})=\sigma_i(M)$ and so
$A_iB_i=T_i$. 
Also, 
\begin{multline*}
A_i\cap B_i=\sigma_i(K_{j_1})\cap\sigma_i(K_{j_2})=\sigma_i\left(\overline 
K_{j_1}\right)\cap\sigma_i\left(\overline K_{j_2}\right)\\=\sigma_i\left(\overline 
K_{1}\right)\cap\cdots\cap\sigma_i\left(\overline K_{\ell}\right)=\sigma_i\left(\overline K_1\cap\cdots\cap \overline K_\ell\right)=\sigma_i\left(\overline M_\omega\right).
\end{multline*}
Hence {\textsf{Prop2\!$\not\sim$[b]}(iii)  is valid. 
If
$g\in\norm{G_\omega}{A_i}$ then $T_i^g\cap T_i\neq 1$ and so $T_i^g=T_i$. 
Thus $g\in\norm{G_\omega}{T_i}$, and so $\norm{G_\omega}{A_i}\leq \norm{G_\omega}{T_i}$. Similarly $\norm{G_\omega}{B_i}\leq \norm{G_\omega}{T_i}$. Suppose now that  $g\in\norm{G_\omega}{T_i}$.  
Then $\sigma_i(
K)^g=\sigma_i(K^g)$ for all $K\in\K$. Thus the uniqueness of $j_1$ and $j_2$ yields that $g$ fixes $K_{j_1}$ and $K_{j_2}$. Therefore $g$
fixes $\sigma_i(K_{j_1})$ and $\sigma_i(K_{j_2})$, and so $g\in\norm{G_\omega}{A_i}\cap\norm{G_\omega}{B_i}$. Thus all
properties in {\sf Prop2\!$\not\sim$[b]} hold.

\textbf{\textsf{Prop2\!$\not\sim$[c]}}\quad
Let $\Gamma$ denote $\Gamma(G,\E)$.
It follows from the definition of $\Gamma$ that
the conjugation action of $G_\omega$ 
is transitive on the vertex-set $\K$ and on the
arc-set $\T$ of
$\Gamma$. We claim that $G_\omega$ preserves adjacency in
$\Gamma(G,\E)$. Let $i_1\in\{1,\ldots,k\}$ 
with $\beta(T_{i_1})=K_{j}$, so that $\sigma_{i_1}(K_j)$ is 
$G_\omega$-conjugate to $A_1$, 
and let $g\in G_\omega$. Let
$i_2\in\{1,\ldots,k\}$ such that $T_{i_1}^g=T_{i_2}$. Then
$\sigma_{i_1}(K_{j})^g=\sigma_{i_2}(K_{j}^g)$.
As
$\sigma_{i_1}(K_j)$ is 
$G_\omega$-conjugate to $A_1$, so is $\sigma_{i_2}(K_j^g)$.
Thus $\beta(T_{i_1}^g)=\beta(T_{i_2})=K_{j}^g=\beta(T_{i_1})^g$,
as required. Thus $\beta$ is preserved by the $G_\omega$-action; similar argument shows that $\varepsilon$ is also preserved by the $G_\omega$-action. Hence
all claims of the theorem hold.\hfill$\Box$

\medskip

\subsection{A converse of Theorem~\ref{main2cnot}}\label{7.1}\quad
Theorem~\ref{main2cnot} can be reversed in the following sense.
Let $G$ be a finite innately transitive group on $\Omega$ with a non-abelian plinth
$M$, and let $T_1,\ldots,T_k$ be the simple direct factors of $M$. Assume that a point stabiliser
$M_\omega$ can be decomposed as
$M_\omega=\sigma_1(M_\omega)\times\cdots\times\sigma_k(M_\omega)$. Set
$\T=\{T_1,\ldots,T_k\}$. 
Suppose that $A_1,\ B_1$ are subgroups of
$T_1$ and $\Gamma=(V,\T,\beta,\varepsilon)$ is a generalised di-graph, such
that properties \textsf{Prop2\!$\not\sim$[b]} and \textsf{Prop2\!$\not\sim$[c]} hold.
This 
amounts to saying that 
\begin{enumerate}
\item[(i)] $A_1,\ B_1$
are proper subgroups of $T_1$;
\item[(ii)] $A_1$ and $B_1$ are
not conjugate under $G_\omega$;
\item[(iii)] $A_1B_1=T_1$,  $A_1\cap B_1=\sigma_1(M_{\omega})$;
\item[(iv)] $\norm {G_\omega}{T_1}=\norm{G_\omega}{A_1}=\norm{G_\omega}{B_1}$;
\end{enumerate}
and also that 
$G_\omega$ induces a vertex and arc-transitive group of automorphisms
of $\Gamma$, where the $G_\omega$-actions are by
conjugation.

For $i=1,\ldots,k$, choose $g_i\in G_\omega$ such that
$T_1^{g_i}=T_i$. 
For each element $v\in V$ set $K_v=\prod_{i=1}^k K_{v,i}$ where 
$$
K_{v,i}=\left\{\begin{array}{ll}
A_1^{g_i} & \mbox{if $\beta(T_i)=v$;}\\
B_1^{g_i} & \mbox{if $\varepsilon(T_i)=v$;}\\
T_i & \mbox{otherwise.}\end{array}\right.
$$
First we prove that the $K_{v,i}$ are well-defined, that is, their definitions
are independent of the choice of the $g_i$. Suppose that $g_i,\ g_i'\in
G_\omega$ are such that $T_1^{g_i}=T_1^{g_i'}=T_i$. Then
$g_ig_i'{}^{-1}\in\norm {G_\omega}{T_1}$ and so, by property~(iv) above,  
$g_ig_i'{}^{-1}\in\norm
{G_\omega}{A_1}\cap\norm{G_\omega}{B_1}$. Hence $A_1^{g_i}=A_1^{g_i'}$ and
$B_1^{g_i}=B_1^{g_i'}$. Thus the $K_{v,i}$ are well-defined.

Let $\K=\{K_v\ |\ v\in V\}$. 
We claim that $\K$ is a $G_\omega$-invariant
Cartesian system for $M$ with respect to $\omega$. 
First note that the $K_v$ are direct products of their projections
and, for all $i$,
$$
\bigcap_{v\in V}K_{v,i}=A_1^{g_i}\cap B_1^{g_i}=(A_1\cap
B_1)^{g_i}=\sigma_1(M_\omega)^{g_i}=\sigma_{i}(M_\omega).
$$
Therefore 
$$
\bigcap_{v\in V}K_v=\prod_{i=1}^k\sigma_i(M_\omega)=M_\omega.
$$
Hence~\eqref{csdef1} holds.
The choice of $A_1$ and $B_1$ and the definition of the $K_{v,i}$ imply that for each $i$ and $v$, 
$$
K_{v,i}\left(\bigcap_{v'\neq v} K_{v',i}\right)=T_i.
$$
As $K_{v,i}\leq K_v$ for all $i$ and $v$, it follows that 
$K_v\left(\bigcap_{v'\neq v}K_{v'}\right)=M$ for all
$v$. Thus \eqref{csdef2}~holds and $\K$ is a Cartesian system for
$M$ with respect to $\omega$.  
We prove now that the set $\K$  is invariant under
conjugation by $G_\omega$. 
Let $v\in V$, $i_1,\ i_2\in\{1,\ldots,k\}$ and $g\in G_\omega$ such
that $T_{i_1}^g=T_{i_2}$. 
We claim that $K_{v,i_1}^g=K_{v^g,i_2}$. Suppose first that
$\beta(T_{i_1})=v$. Then, as $g$ induces
an automorphism of $\Gamma$, we obtain that
$\beta(T_{i_1}^g)=v^g$, that is, $\beta(T_{i_2})=v^g$. Thus in this
case we have $K_{v^g,i_2}=A_1^{g_{i_2}}$,
$K_{v,i_1}=A_1^{g_{i_1}}$, and hence
$K_{v,i_1}^g=A_1^{g_{i_1}g}$. As 
$T_1^{g_{i_1}g}=T_{i_1}^g=T_{i_2}$, we obtain, as above, that $A_1^{g_{i_1}g}
=A_1^{g_{i_2}}$. Hence $K_{v,i_1}^g=K_{v^g,i_2}$. Similar argument shows that 
this equality also
holds when $\varepsilon(T_{i_1})=v$, and when $T_{i_1}$ is not
adjacent to $v$. 
Therefore 
$$
K_v^g=\left(\prod_{i=1}^k K_{v,i}\right)^g=\prod_{i=1}^k K_{v^g,i}=K_{v^g}.
$$
Hence $\K$ is $G_\omega$-invariant. The last equation also
shows that the $G_\omega$-actions on $V$ and on $\K$ are equivalent. Thus $\K$ is a
$G_\omega$-transitive Cartesian system of subgroups for $M$ with
respect to $\omega$, and it follows from the definition of the $K_v$
that $\E(\K)\in\cd {2\not\sim}G$, $\Gamma=\Gamma(G,\E)$, and 
$\mathcal F_i(\E(\K),M,\omega)=\{A_i,B_i\}$ for each $i$. 

\medskip

As in the previous section, it is possible to investigate further
the conditions that ensure the existence of a Cartesian decomposition
in $\cd{2\not\sim}G$ for some innately transitive group $G$. If this set
of decompositions is non-empty then the
set of simple direct factors of the plinth must play the r\^ole of the
arc-set of a generalised di-graph. Thus we expect that the nature of
the conjugation action
on the simple direct factors can strongly restrict the possible generalised 
di-graphs satisfying \textsf{Prop2\!$\not\sim$[c]}, and hence the possible elements of $\cd{2\not\sim}G$. 
Though details of this phenomenon
are not addressed in this paper, we present a simple example
for illustration.

\begin{example}
We claim that no generalised di-graph exists having four~arcs and admitting 
an automorphism group that
acts vertex and arc-transitively inducing an $\alt 4$ or $\sy 4$ on 
the arcs. For, if 
$\Gamma$ is such a generalised di-graph then every vertex has a constant 
number of outgoing arcs. Hence the number of vertices must be a divisor of~4
(and is not 1 by the definition of a generalised di-graph).
It is left to the reader to check that no such graph exists on 2~or  
4~vertices. Therefore if 
$G$ is a finite  innately transitive group with a non-abelian plinth 
$M=T_1\times\cdots\times T_4$ where the permutation action of $G$ on 
the $T_i$ is permutationally isomorphic to $\alt 4$ or $\sy 4$ then 
$\cd{2\not\sim}G=\emptyset$.
\end{example}

\section{Cartesian decompositions in $\cd{{\rm 1S}}{G}$}\label{6}

In this section the following notation is used. 
Let $G$ be a finite innately transitive group on $\Omega$ 
with a non-abelian plinth $M$,
and let $T_1,\ldots,T_k$ be the simple normal subgroups of $M$, each
isomorphic to the simple group $T$. Let
$\omega\in\Omega$, and set $\T=\{T_1,\ldots,T_k\}$.  
Let $\overline M_\omega=\norm{M}{M_\omega}$ and let $\overline\Omega$ 
denote the $M$-invariant partition $\parti M{\overline M_\omega}$ of $\Omega$. 
Then $\overline M_\omega$ is normalised by $G_\omega$. 
Thus Lemma~\ref{quotlem} shows that $\overline\Omega$ is $G$-invariant.
Let $\overline\omega$ be the block in
$\overline\Omega$ that contains $\omega$. Then Lemma~\ref{quotlem}
also implies that $G_{\overline\omega}=\overline M_\omega G_\omega$, that
$\overline\omega$ is the $\overline M_\omega$-orbit containing $\omega$,
and that $\overline M_\omega = M_{\overline\omega}$. Let $\overline G$ denote the 
permutation group on $\overline\Omega$ induced by $G$, so $\overline 
G_{\overline\omega}$ is the subgroup of $\overline G$ induced by $G_{\overline
\omega}$.

Suppose that $\E\in\cd{1S}G$ and let $\K_\omega(\E)=\{K_1,\ldots,K_\ell\}$ be
the corresponding Cartesian system. For $K\in\K_\omega(\E)$, let $\X_K$ denote 
the set of non-trivial, full strips involved in $K$, and set $\X=\X_{K_1}\cup
\cdots\cup\X_{K_\ell}$. By Theorem~\ref{stripth}, 
$\X$ contains $k/2$ pairwise disjoint, full strips, each of length~2.
Let
$$
\overline K_i=\prod_{X\in\X_{K_i}} X\times\prod_{T_m\not\in\bigcup_{X\in\X_{K_i}}\supp
  X}\sigma_m(K_i).
$$
Set $\overline \K_\omega(\E)=\{\overline K_1,\ldots,\overline K_\ell\}$. 
If $X$ is a strip in $M$ then we define $\min
X=\min\{i\ |\ T_i\in \supp X\}$ and $\max X=\max\{ i \ |\ T_i\in \supp X\}$. 
Suppose that $\X=\{X_1,\ldots,X_{k/2}\}$, and, for $i=1,\ldots,k/2$, let 
$A_i$ and $B_i$ be defined as follows. 
There are unique indices $j_1$ and $j_2$ such that $\sigma_{\min X_i}
(K_{j_1})\neq T_{\min X_i}$ and $\sigma_{\max X_i}(K_{j_2})\neq T_{\max X_i}$; 
set $A_i=\sigma_{\min X_i}(K_{j_1})$ and $B_i=\sigma_{\max X_i}(K_{j_2})$.
Let $\Gamma(G,\E)$ denote the graph $(\K_\omega(\E)\cup\X,E_1\cup E_2)$ where, 
for $K\in\K_\omega(\E)$ and $X\in\X$, $\{K,X\}\in E_1$ if either 
$\sigma_{\min X}(K)< T_{\min X}$ or
$\sigma_{\max X}(K)<T_{\max X}$, and $\{K,X\}\in E_2$ if
$X$ is involved in $K$.

\begin{theorem}\label{main1s}
Let $G$ and $M$ be as in the first paragraph of this section. If $\E\in\cd {1S}G$, then the properties
{\sf Prop1S[a]}--{\sf Prop1S[d]} below all hold.  
\end{theorem}

\noindent\textbf{\textsf{Prop1S[a]} (Quotient Action Property).}
The group $M$ is faithful on $\overline\Omega$, and so, if $K$ is a subgroup of $M$, then we
identify $K$ with its image under the action on $\overline\Omega$. 
The set $\overline\K_\omega(\E)$ is a $\overline G_{\overline\omega}$-invariant
Cartesian system of subgroups for $M$ with respect to
$\overline\omega$. Moreover, $\E(\overline\K_\omega(\E))\in\cd {1S}{\overline G}$. 

\bigskip

\noindent\textbf{\textsf{Prop1S[b]} (Factorisation Property).}
If $i\in\{1,\ldots,k/2\}$ then
\begin{enumerate}
\item[(i)] $X_i$ is a full strip of length~2;
\item[(ii)] $A_i$ is a proper subgroup of $T_{\min X_i}$ and 
$B_i$ is a proper subgroup of $T_{\max X_i}$; 
\item[(iii)] $A_i$ and $B_i$ are
conjugate under $G_\omega$;
\item[(iv)] $X_i(A_i\times B_i)=T_{\min X_i}\times T_{\max X_i}$,   $X_i\cap
(A_i\times B_i)=\sigma_{\supp X}(\norm M{M_{\omega}})$;
\item[(v)] $\norm {G_\omega}{T_{\min X_i}\times T_{\max X_i}}=\norm{G_\omega}
{X_i}=\norm{G_\omega}{A_i\times B_i}$.
\end{enumerate}

\bigskip

\noindent\textbf{\textsf{Prop1S[c]} (Combinatorial Property).}
The graph $\Gamma(G,\E)$ is bipartite, with bipartition formed
by the sets $\K_\omega(\E)$ and $\X$. 
The group $G_\omega$ induces a group of automorphisms of
the graph $\Gamma(G,\E)$, such that $\K_\omega(\E)$, $\X$, $E_1$, 
and $E_2$ are $G$-orbits. 
Moreover, each element of $\X$ is adjacent to one edge 
or two edges from $E_1$, and
 one edge of $E_2$.
Further, if the elements of $\X$ are adjacent to two elements of
$E_1$ then the following condition must also hold: if, for some $K\in\K_\omega(\E)$ and $X\in\X$,  $\{K,X\}\in E_1$
then $(G_\omega)_{\{K,X\}}=\norm{G_\omega}{T_{\min X}}\cap \norm{G_\omega}{T_{\max
X}}$. 

\bigskip

\noindent\textbf{\textsf{Prop1S[d]} (Isomorphism Property).}
The group $T$, the subgroup $A_i$ (which is $G_\omega$-conjugate to $B_i$), 
and $\sigma_{\supp X_i}(M_{\overline\omega})$ are as in
Table~\ref{isomfact}. 
The group $\overline G$ is permutationally isomorphic to a subgroup of
$\aut M$ acting on $\overline\Omega$. In particular, $M$ is the unique minimal
normal subgroup of $\overline G$, and  $\overline G$ is
quasiprimitive. Moreover, if $T$ is as in rows
1--3 of Table~\ref{isomfact} then $\overline M_\omega=M_\omega$, 
and  $G\cong\overline G$ as permutation groups. Otherwise each block 
in $\overline\Omega$ has size dividing $2^{k/2}$, the kernel $N$ of
the action of $G$ on $\overline\Omega$ is an elementary abelian
2-group of rank at most $k/2$ and $\overline G\cong G/N$. 

\bigskip

\subsection{Some computation}
The following three results are needed for the proof of Theorem~\ref{main1s}.
The results are stated and proved in the context 
introduced in the first two paragraphs of this section.

\begin{proposition}\label{6.1th}
If $\E\in\cd{\rm 1S}G$,
then the group $T$ is as in one of the rows of
Table~$\ref{isomfact}$, and each $\mathcal F_i(\E,M,\omega)$
contains a subgroup isomorphic to the group $A_i$ in the
corresponding row of Table~$\ref{isomfact}$. 
Further, for all $K\in\K_\omega(\E)$, 
\begin{equation}\label{weak2}
\prod_{X\in\X_K}X\times\prod_{T_m\not\in\bigcup_{X\in\X_K}\supp
X}\sigma_m(K)'\ \leq\ K.
\end{equation}
If  $T$ is     as in rows~$1$--$3$ then  for all $i$,
\begin{equation}\label{strong2}
K=\prod_{X\in\X_K}X\times\prod_{T_m\not\in\bigcup_{X\in\X_K}\supp
X}\sigma_m(K)\quad\mbox{for all}\quad K\in\K_\omega(\E).
\end{equation}
\end{proposition}

\begin{center}
\begin{table}[ht]
$$
\begin{array}{|l|c|c|c|}
\hline
 & T & A_i &  \sigma_{\supp X_i}(M_{\overline\omega})\\
\hline
1 & \alt 6&\alt 5 & \dih{10}\\
\hline
2 & \mat{12}&\mat{11} & \psl 2{11}
 \\
\hline
3 & \pomegap 8q &\Omega_7(q) & {\sf G}_2(q)\\
\hline
4 & \sp 4q,\ q\geq 4\mbox{ even}
& \sp 2{q^2}\cdot 2 & \dih{q^2+1}\cdot 2 \\
\hline
\end{array}
$$
\caption{Factorisations of finite simple groups with two isomorphic subgroups}
\label{isomfact}
\end{table}
\end{center}

\begin{proof}
Let $\K_\omega(\E)=\{K_1,\ldots,K_\ell\}$, and for $i=1,\ldots,\ell$
let $\widehat K_{i}$ denote
$\bigcap_{j\neq i}K_j$. 
By Theorem~\ref{stripth}, each non-trivial, full strip involved in a
$K_i$ has length~2. Suppose without loss of generality that $X$ is a
non-trivial full strip involved in $K_1$ covering
$T_1$ and $T_2$. 
Thus by (\ref{csdef2}),
$T_{1}\times T_{2}=\sigma_{\{1,2\}}(K_1)\sigma_{\{1,2\}}(\widehat K_{{1}})$. 
Suppose that $X=\{(t,\alpha(t))\
|\ t\in T_1\}$ for some isomorphism $\alpha:T_1\rightarrow T_2$.  
Then it follows from~\cite[Lemma~2.1]{charfact} that
$T_1=\sigma_1(\widehat K_{{1}})\alpha^{-1}(\sigma_2(\widehat
K_{{1}}))$. By the definition of $\cd {1S}G$, $\sigma_1(K_{j_1})\neq T_1$ and
$\sigma_2(K_{j_2})\neq T_2$ for some $j_1,\
j_2\in\{2,\ldots,\ell\}$. Thus $\sigma_1(\widehat K_1)$
and $\sigma_2(\widehat K_1)$ are proper subgroups of $T_1$ and $T_2$,
respectively. 
Moreover, if $g\in G_\omega$ such that $T_1^g=T_2$ then $T_2\in\supp X\cap
\supp X^g$, and so Theorem~\ref{stripth} implies that
$X^g=X$. Hence, again by Theorem~\ref{stripth}, $g\in\norm{G_\omega}{K_1}$, 
and also $g\in\norm{G_\omega}{\widehat K_{1}}$. Thus $\sigma_1(\widehat
K_{1})^g=\sigma_2(\widehat K_{1})$ is a proper subgroup of $T_2$. 
Hence $T_1=\sigma_1(\widehat
K_{{1}})\alpha^{-1}(\sigma_2(\widehat K_{{1}}))$
is a factorisation with proper, isomorphic
subgroups. Therefore~\cite[Lemma~5.2]{recog} implies that $T_1\cong
T$ is as in Table~\ref{isomfact} and the isomorphism types of
$\sigma_1(\widehat K_{{1}})$ and $\sigma_2(\widehat K_{{1}})$ are as
in the $A_i$-column of the corresponding row of Table~\ref{isomfact}.

Suppose that $\sigma_i(K_j)\neq T_i$
for some $i$ and $j$. Then there is a non-trivial full 
strip $X\in\X$ covering $T_i$; assume that $X\in\X_{K_m}$ for some
$m\in\{1,\ldots,\ell\}\setminus\{j\}$ and that $\supp
X=\{T_i,T_{i'}\}$. Then the argument of
the previous paragraph shows that $\sigma_i(\widehat K_{{m}})$ is as in the $A_i$-column of Table~\ref{isomfact}. In particular,
$\sigma_i(\widehat K_{{m}})$ is a maximal subgroup of $T_i$. Since
$\widehat K_{{m}}\leq
K_j$, we obtain that $\sigma_i(\widehat K_{{m}})\leq\sigma_i(K_j)<T_i$, and so
$\sigma_i(\widehat K_{{m}})=\sigma_i(K_j)$. Therefore $\sigma_i(K_j)$
is also as in the $A_i$-column of the corresponding row of the table.

We have proved so far that for
all $i$ and $j$ 
either $\sigma_i(K_j)=T_i$ or $\sigma_i(K_j)\cong A$ where $A$ is as in
the $A_i$-column of Table~\ref{isomfact}. In particular $A$ is a maximal 
subgroup of $T$, $A$ is
almost simple,  and if $T$ is as in rows~1--3 of Table~\ref{isomfact} then
$A$ is simple. Suppose by contradiction that~\eqref{weak2} fails to
hold for some $K_j$. Set 
$$
\overline S=\left\{m\ |\ T_m\in\bigcup_{X\in\X_{K_j}}\supp X\right\}\quad 
\mbox{and}\quad S=\{1,\ldots,k\}\setminus \overline S,
$$
and write $\sigma_S, \sigma_{\overline S}$ for the projection of $M$ onto 
$\prod_{s\in S}T_s$ and $\prod_{s\in\overline S}T_s$ respectively.
Then it follows from the definition of $\X_{K_j}$ that
$$
K_j=\sigma_{S}(K_j)\times\sigma_{\overline S}(K_j).
$$
As~\eqref{weak2} fails for $K_j$ we must have that
$$
\prod_{m\in S}\sigma_m(K_j)'\not\leq \sigma_{S}(K_j).
$$
Thus it follows from~\cite[Lemma~2.3]{charfact} 
that there are distinct elements $i_1,\ i_2$ of $S$ such
that 
\begin{equation}\label{wrong2}
\sigma_{i_1}(K_j)'\times\sigma_{i_2}(K_j)'\not\leq
\sigma_{\{i_1,i_2\}}(K_j).
\end{equation}
If $\sigma_{i_1}(K_j)=T_{i_1}$ then, by~\cite[Lemma~4.3]{transcs},
there is a full strip $X$ involved in $K_j$ covering $T_{i_1}$. By the
definition of $S$, we must have that $X=T_{i_1}$, and so
$\sigma_{i_1}(K_j)\leq K_j$. Hence
$\sigma_{i_1}(K_j)\leq\sigma_{\{i_1,i_2\}}(K_j)$, and also
$\sigma_{i_2}(K_j)\leq\sigma_{\{i_1,i_2\}}(K_j)$. Hence
$\sigma_{i_1}(K_j)\times \sigma_{i_2}(K_j)=\sigma_{\{i_1,i_2\}}(K_j)$,
contradicting~\eqref{wrong2}. Thus $\sigma_{i_1}(K_j)$ is a proper
subgroup of $T_{i_1}$, and also $\sigma_{i_2}(K_j)$ is a proper
subgroup of $T_{i_2}$.

By Theorem~\ref{stripth}(d), 
$G_\omega$ is transitive on $\X$, and so there are (not necessarily distinct) 
strips $X_1$ and $X_2$ in $\X$ such that $X_1$ covers $T_{i_1}$ and 
$X_2$ covers $T_{i_2}$. 
Suppose that $X_1=X_2$. Then $\supp X_1=\{T_{i_1},T_{i_2}\}$, and
let ${j_1}\in\{1,\ldots,\ell\}\setminus\{j\}$ be such that
$X_1\in\X_{K_{j_1}}$. Then, as verified above, $\sigma_{i_1}(\widehat
K_{j_1})$ and $\sigma_{i_2}(\widehat K_{j_1})$ are maximal subgroups of
$T_{i_1}$ and $T_{i_2}$, respectively, and, in addition, $\sigma_{i_1}(\widehat
K_{j_1})\cong \sigma_{i_2}(\widehat K_{j_1})$. Thus the factorisation
$$
X_1\sigma_{\{i_1,i_2\}}(\widehat K_{j_1})=\sigma_{\{i_1,i_2\}}(K_{j_1})
\sigma_{\{i_1,i_2\}}(\widehat K_{j_1})=T_{i_1}\times T_{i_2}
$$ 
is as
in~\cite[Theorem~1.5]{charfact}. Hence~\cite[Theorem~1.5]{charfact}
implies that 
$$
\sigma_{i_1}(\widehat K_{j_1})'\times\sigma_{i_2}(\widehat K_{j_1})'\leq\sigma_{\{i_1,i_2\}}(\widehat
K_{j_1}).
$$
Note that $j\neq j_1$, and so $\sigma_{\{i_1,i_2\}}(\widehat
K_{j_1})\leq \sigma_{\{i_1,i_2\}}(K_{j})$. Moreover, $\sigma_{i_1}(\widehat K_{j_1})$ is a maximal subgroup of $T_{i_1}$ and so is
$\sigma_{i_1}(K_j)$. As $\sigma_{i_1}(
\widehat K_{j_1})\leq \sigma_{i_1}(K_j)$, we obtain that $\sigma_{i_1}(
\widehat K_{j_1})= \sigma_{i_1}(K_j)$, and, similarly, 
$\sigma_{i_2}(K_{j_1})=\sigma_{i_2}(K_j)$.  
Therefore 
$$
\sigma_{i_1}(K_{j})'\times\sigma_{i_2}(K_{j})'\leq\sigma_{\{i_1,i_2\}}(K_{j}),
$$
which is a contradiction. Hence $X_1\ne X_2$. 

Suppose that $X_1$ is involved in $K_{j_1}$ and $X_2$ is
involved in $K_{j_2}$, where $j_1$ and $j_2$ are not necessarily distinct
elements of $\{1,\ldots,\ell\}\setminus\{j\}$. Let $I=\supp X_1\cup
\supp X_2$ and set $\widehat K_{{j_1,j_2}}=\bigcap_{m\neq j_1,\
j_2}K_m$.  Then, by~\cite[Lemma~3.1]{recog}, $(K_{j_1}\cap
K_{j_2})\widehat K_{j_1,j_2}=M$, and so
$$
\sigma_I(M)=\sigma_I(K_{j_1}\cap
K_{j_2})\sigma_I\left(\widehat K_{j_1,j_2}\right).
$$

Suppose
that $n\in\supp X_1\cup \supp X_2$; in fact suppose without loss of
generality that $n\in\supp X_1$. Then the argument above shows that
$\sigma_n(\widehat K_{{j_1}})\cong A$ and also $\sigma_n(K_{j'})\cong
A$ where $A$ is as in the $A_i$-column of Table~\ref{isomfact} and 
$j'\in\{1,\ldots,\ell\}$ is such that $\sigma_n(K_{j'})<T_n$. Since, 
$$
\sigma_n(\widehat K_{{j_1}})\leq\sigma_n(\widehat K_{{j_1,j_2}})\leq\sigma_n(K_{j'}),
$$
we obtain that $\sigma_n(\widehat K_{{j_1,j_2}})\cong A$, and this holds  
for all $n\in\supp X_1\cup\supp X_2$. 
Clearly $\sigma_I(K_{j_1}\cap K_{j_2})\leq$ $ X_1\times X_2$,
and so 
$$
\sigma_I(M)=(X_1\times X_2)\sigma_I\left(\widehat K_{{j_1,j_2}}\right).
$$
Then it follows from~\cite[Theorem~1.5]{charfact} that
$$
\sigma_{\min X_1}(\widehat K_{{j_1,j_2}})'\times \sigma_{\max
  X_1}(\widehat K_{{j_1,j_2}})'\times \sigma_{\min
  X_2}(\widehat K_{{j_1,j_2}})'\times \sigma_{\max
  X_2}(\widehat K_{{j_1,j_2}})'\leq\sigma_{I}(\widehat K_{{j_1,j_2}}).
$$
As $i_1,\ i_2\in I$, we obtain that
$$
\sigma_{i_1}(\widehat K_{{j_1,j_2}})'\times
\sigma_{i_2}(\widehat K_{{j_1,j_2}})'\leq
\sigma_{\{i_1,i_2\}}(\widehat K_{{j_1,j_2}})\leq \sigma_{\{i_1,i_2\}}(K_j).
$$
Since $\sigma_{i_1}(\widehat K_{{j_1,j_2}})'=\sigma_{i_1}(K_j)'$ and
$\sigma_{i_2}(\widehat K_{{j_1,j_2}})'=\sigma_{i_2}(K_j)'$, this is a contradiction.
Hence~\eqref{weak2} holds. 
If $T$ is as in rows~1--3, then
$\sigma_i(K_j)$ is simple, and hence perfect. This
proves~\eqref{strong2}.
\end{proof}

Next we need to compute normalisers of point stabilisers and 
Cartesian system elements.

\begin{lemma}\label{knorm}
If $\E\in\cd{1S}G$, then for all $i\in\{1,\ldots,\ell\}$, 
$$
\norm
M{K_i}=\overline K_i\quad\mbox{and}\quad
K_i'=\prod_{X\in\X_{K_i}}X\times\prod_{T_j\not\in\bigcup_{X\in\X_{K_i}}\supp X}\sigma_j(K_i)'.
$$
\end{lemma}
\begin{proof}
If $T$ is as in rows~1--3 of Table~\ref{isomfact} then the claim of
the lemma follows from the fact that, by Lemma~\ref{normstriplemma}, 
each strip $X\in\X_{K_i}$ is
a simple and self-normalising subgroup of $\sigma_{\supp X}(M)$, and if $T_j$ is not covered
by any strip in $\X_{K_i}$ then $\sigma_j(K_i)$ is self-normalising in
$T_j$. 

Suppose now that $T$ is as in row~4 of Table~\ref{isomfact} and set
$$
\underline K_i=\prod_{X\in\X_{K_i}} X\times\prod_{T_j\not\in\bigcup_{X\in\X_{K_i}}\supp
  X}\sigma_j(K_i)'.
$$ 
Then it follows from Proposition~\ref{6.1th} that $\underline K_i\leq K_i$,
and from the definition of $\overline K_i$ that $K_i\leq \overline
K_i$. Now each strip $X\in\mathcal X_{K_i}$ is self-normalising in $\sigma_{\supp X}(M)$, 
and $\norm {T_j}{\sigma_j(K_i)'}=\norm{T_j}{\sigma_j(K_i)}=
\sigma_j(K_i)$ whenever
$T_j\not\in\supp X$ for some $X\in\X_{K_i}$. Hence $\norm M{\underline
  K_i}=\overline K_i$ and $\norm M{K_i}\leq\overline K_i$.
On the other hand, Lemma~\ref{normlemma1}
implies that $\norm M{\underline K_i}\leq\norm M{K_i}$, and so 
$\overline K_i = \norm M{\underline K_i}=\norm M{K_i}$. 

It remains to prove that $K_i'=\underline K_i$. As for $j=1,\ldots,k$, 
either $\sigma_j(K_i)'=\sigma_j(K_i)$ or $\sigma_j(K_i)/\sigma_j(K_i)'$ is
isomorphic to $\Z_2$, it follows that
$K_i/\underline K_i$ is an elementary abelian $2$-group. Thus $K_i'\leq
\underline K_i$. On the other hand, $\underline K_i$ is a direct product
of non-abelian, finite simple groups, and so no quotient of
$\underline K_i$ is abelian. This
proves that $K_i'=\underline K_i$ as required.
\end{proof}

Suppose that $G_1,\ G_2$ are groups and let $H$ be a subgroup of
$G_1$. If $\alpha:H\rightarrow G_2$ is an injective homomorphism then
we define 
$$
\diag\alpha=\{(h,\alpha(h))\ |\ h\in H\}
$$
as a subgroup of $G_1\times G_2$.

\begin{lemma}\label{uniqueext}
Let $\E\in\cd{\rm 1S}G$ and let $\K_\omega(\E)=\{K_1,\ldots,K_\ell\}$. Then 
$$
\norm M{M_\omega}=\norm M{K_1}\cap\cdots\cap\norm M{K_\ell},
$$
and
$\norm M{M_\omega}$ is the direct product
of $k/2$ strips of length~$2$ in $M$. Moreover, if $Y$ is a strip in
$\norm M{M_\omega}$ then $Y$ is isomorphic to  
the group in the last column of the appropriate row of Table~$\ref{isomfact}$.
\end{lemma}
\begin{proof}
For $i=1,\ldots,\ell$ we
have $K_i'\leq K_i\leq\norm M{K_i}$ with 
equality if $T$ is as in one of the rows~1--3 of Table~\ref{isomfact} 
(see Lemma~\ref{knorm}). Let
$X_1,\ldots,X_{k/2}$ be the strips involved in $\K_\omega(\E)$. Then
Theorem~\ref{stripth} implies that $\{\supp X_1,\ldots,\supp
X_{k/2}\}$ is a partition of $\{T_1,\ldots,T_k\}$, and it follows from
Lemma~\ref{knorm} that, for $i=1,\ldots,\ell$,
$$
K_i'=\prod_{j=1}^{k/2}\sigma_{\supp X_j}(K_i')\quad\mbox{and}\quad
\norm M{K_i}=\prod_{j=1}^{k/2}\sigma_{\supp X_j}(\norm M{K_i}).
$$
Let $\underline{M}_\omega=K_1'\cap\cdots\cap K_\ell'$. Set 
$M_0=\overline K_1\cap\cdots\cap\overline K_\ell$ and recall that 
$M_0=\norm M{K_1}\cap\cdots\cap\norm M{K_\ell}$. Then 
$$
\underline M_\omega=\prod_{j=1}^{k/2}\sigma_{\supp X_j}(\underline M_\omega)
\quad\mbox{and}\quad
M_0=\prod_{j=1}^{k/2}\sigma_{\supp X_j}(M_0).
$$ 
For $i=1,\ldots,k/2$, the subgroup $X_i=\diag\alpha_i$ for some isomorphism 
$\alpha_i:T_{\min X_i}\rightarrow T_{\max X_i}$. 
Let $K_{j_1}$ and $K_{j_2}$ be the elements of the Cartesian system
such that
$\sigma_{\min X_i}(K_{j_1})\neq T_{\min X_i}$ and $\sigma_{\max
X_i}(K_{j_2})\neq T_{\max X_i}$. 
Set 
$$
\widehat Y_i=\sigma_{\min X_i}(K_{j_1})\cap\alpha_i^{-1}(\sigma_{\max
X_i}(K_{j_2}))\quad\mbox{and}\quad \check Y_i=\sigma_{\min X_i}(K_{j_1})'\cap\alpha_i^{-1}(\sigma_{\max
X_i}(K_{j_2}))'. 
$$
Let $\widehat \alpha_i$ and
$\check\alpha_i$ denote
the restrictions of $\alpha_i$ to the subgroups $\widehat Y_i$ and $\check
Y_i$, respectively.
Then we have that 
$$
\sigma_{\supp
X_i}(M_0)=\diag{\widehat \alpha_i}\quad\mbox{and}\quad
\sigma_{\supp
X_i}(\underline M_\omega)=\diag{\check\alpha_i}.
$$

Suppose first that $T$ is as in rows~1--3 in
Table~\ref{isomfact}. Then, as $M_\omega=\underline M_\omega=M_0$, it follows that $M_\omega$ is the direct product of the
$\diag{\widehat \alpha_i}$. 
On the other hand, by~\cite[Lemma~2.1]{charfact}, the factorisation 
$T_{\min X_i}=\sigma_{\min X_i}(K_{j_1})\alpha_i^{-1}(\sigma_{\max
X_i}(K_{j_2}))$ involves isomorphic subgroups. 
Thus the subgroups involved in
this factorisation must be as
in~\cite[Lemma~5.2]{recog}. Now the isomorphism type of the
intersection $\widehat Y_i$
can be determined using the~\cite{atlas} in rows 1--2 and
\cite[3.1.1(vi)]{kleidman} in row~3. Hence we
find that, for $T$ in one of these rows, the group $\widehat Y_i$ is 
isomorphic to the subgroup in the last column of
Table~\ref{isomfact}. By~\cite[Lemma~5.2]{recog}, $\widehat Y_i$ is self-normalising with trivial centraliser in
$T_{\min X_i}$, and so Lemma~\ref{normstriplemma} implies that $\norm
{T_{\min X_i}\times T_{\max
X_i}}{\diag\widehat\alpha}=\diag\widehat\alpha$, and so $\norm
M{M_\omega}=M_\omega$, as required.

Suppose now that $T$ is as in row~4 of Table~\ref{isomfact}. Then the 
isomorphism $\widehat Y_i\cong\dih{q^2+1}\cdot 2$ and $\check Y_i\cong
\dih{q^2+1}$ follow
from \cite[3.2.1(d)]{lps:max}. 
Using Lemma~\ref{simpnorm} we obtain that
$\norm{T_{\min X_i}}{\check Y_i}=\norm{T_{\min X_i}}{\widehat Y_i}=\widehat
Y_i$ and $\cent{T_{\min X_i}}{\check Y_i}=\cent{T_{\min X_i}}{\widehat Y_i}=1$.
Thus Lemma~\ref{normstriplemma} implies that
$\norm M{\underline M_\omega}=M_0$, $\underline
M_\omega\cong (\dih{q^2+1})^{k/2}$ and $M_0\cong(\dih{q^2+1}\cdot 2)^{k/2}$, and $M_0/\underline M_\omega$ is an elementary abelian
2-group. Hence $\norm M{M_\omega}\leq M_0$. 
On the other hand
Lemma~\ref{normlemma1} implies that
$ {M_0}\leq\norm M{M_\omega}$. Therefore  $\norm
M{M_\omega}=M_0$, as required. 
\end{proof}

Now we can prove Theorem~\ref{main1s}. 

\subsection{Proof of Theorem~\ref{main1s}}
\textbf{\textsf{Prop1S[a]}}\quad
For each $i$ there is a unique $j$ such that $\sigma_i(K_j)<T_i$. 
Thus the $\overline K_j$ are proper subgroups of $M$, and no $T_i$ is contained in $\overline M_\omega$. 
Hence $M$ acts faithfully on $\overline\Omega$. 
Lemmas~\ref{knorm} and~\ref{uniqueext} imply that
$\overline K_1\cap\cdots\cap\overline K_\ell=\norm M{M_\omega}=\overline 
M_{\omega}$.
Therefore~\eqref{csdef1} holds for $\overline\K_\omega(\E)$. 
As~\eqref{csdef2}
holds for $\K_\omega(\E)$, and, for $i=1,\ldots,\ell$, $K_i\leq\overline K_i$, we
have that~\eqref{csdef2} also holds for $\overline\K_\omega(\E)$. Therefore $\overline\K_\omega(\E)$
is a Cartesian system of subgroups for $M$ with respect to
$\overline\omega$. We claim that $\overline\K_\omega(\E)$ is invariant under conjugation by
$\overline G_{\overline\omega}$. 
Note that $G_{\overline\omega}=M_{\overline\omega}G_\omega$, and so it suffices
to prove that $\overline\K_\omega(\E)$ is invariant under conjugation by
$G_\omega$. This however follows from the fact that
$\{K_1,\ldots,K_\ell\}$ is $G_\omega$-invariant and, by Lemma~\ref{knorm}, 
$\overline\K_\omega(\E)=\{\norm
M{K_1},\ldots,\norm M{K_\ell}\}$. 

\textbf{\textsf{Prop1S[b]}}\quad
It follows from Theorem~\ref{stripth} that
{\sf Prop1S[b]}(i)
holds. It is clear that {\sf Prop1S[b]}(ii)
also holds. Recall that $\overline M_\omega=
M_{\overline\omega}=\norm M{M_\omega}$. Let $X$ be a
non-trivial, full strip involved in $K_i$,  say, and let $j_1,\
j_2\in\{1,\ldots,\ell\}\setminus\{i\}$ be such that $\sigma_{\min
X}(K_{j_1})<T_{\min X}$ and $\sigma_{\max X}(K_{j_2})<T_{\max
X}$. Set $A=\sigma_{\min
X}(K_{j_1})$ and $B=\sigma_{\max X}(K_{j_2})$.
Suppose that $g\in G_\omega$ is such that $T_{\min X}^g=T_{\max X}$. Then 
$A^g=\sigma_{\min X}(K_{j_1})^g=\sigma_{\max X}(K_{j_1}^g)$. As $j_2$ is the 
unique integer such that $\sigma_{\max X}(K_{j_2})<T_{\max
X}$, we obtain that $K_{j_1}^g=K_{j_2}$, and so $A^g=B$. Hence \textsf{Prop1S[b]}(iii) also holds.
Note that, as $\K_\omega(\E)$ is a
Cartesian system, we have that $K_i(K_{j_1}\cap K_{j_2})=M$. Thus 
$$
T_{\min X}\times T_{\max X}=\sigma_{\supp X}(K_j)\sigma_{\supp
X}(K_{j_1}\cap K_{j_2}).
$$
As $\sigma_{\supp
X}(K_{j_1}\cap K_{j_2})\leq A\times B$ we obtain that $T_{\min
X}\times T_{\max X}=X(A\times B)$. Since each $\overline K_j$ is the
direct product, over $X_i\in\X$, of its projection under
$\sigma_{\supp X_i}$, so is the subgroup $\overline M_\omega$. 
Thus $X\cap(A\times B)=\sigma_{\supp X}(\overline M_\omega)=\sigma_{\supp X}(\norm M{M_\omega})$, by Lemma~\ref{uniqueext}. Therefore
{\sf Prop1S[b]}(iv) holds. Let
us now prove {\sf Prop1S[b]}(v). As $A\times B\leq T_{\min X}\times
T_{\max X}$,  $X\leq T_{\min X}\times
T_{\max X}$, and $\{T_{\min X},T_{\max X}\}$ is a block for
the $G_\omega$-action on $\T$, it follows that
$$
\norm{G_\omega}{A\times B}\leq \norm{G_\omega}{T_{\min X}\times
T_{\max X}}\quad\mbox{and}\quad\norm{G_\omega}{X}\leq \norm{G_\omega}{T_{\min X}\times
T_{\max X}}.
$$ 
Let $g\in \norm{G_\omega}{T_{\min X}\times
T_{\max X}}$. Then $X^g$ is a strip involved in  $K_i^g\in\K_\omega(\E)$ such
that $X$ and $X^g$ have the same support. Hence Theorem~\ref{stripth} implies
that $X=X^g$, and so $g\in\norm{G_\omega}X$. Also the element $g$
either normalises both subgroups $T_{\min X}$ and $T_{\max X}$ or swaps these
two subgroups. Hence one of the following scenario holds: either 
$$\sigma_{\min
X}(K_{j_1})^g=\sigma_{\min X}(K_{j_1}^g)\quad\mbox{and} \quad \sigma_{\max
X}(K_{j_2})^g=\sigma_{\max X}(K_{j_2}^g);$$
 or 
$$\sigma_{\min
X}(K_{j_1})^g=\sigma_{\max X}(K_{j_1}^g)\quad\mbox{and}\quad\sigma_{\max
X}(K_{j_2})^g=\sigma_{\min X}(K_{j_2}^g). 
$$
Since $K_{j_1}$ and
$K_{j_2}$ are the unique elements of $\K_\omega(\E)$ whose projection
to $T_{\min X}$ and $T_{\max X}$, respectively, are proper, we obtain
that $\{A^g,B^g\}=\{A,B\}$. Therefore $(A\times B)^g=A\times B$, and so
$g\in\norm{G_\omega}{A\times B}$.

\textbf{\textsf{Prop1S[c]}}\quad
First we prove that $G_\omega$ induces a group of automorphisms of the
graph $\Gamma(G,\E)$. Suppose that, for some $K\in\K_\omega(\E)$ and $X\in\X$, the edge $\{K,X\}$ is in $E_1$ and $g\in
G_\omega$. Then $\sigma_{\min X}(K)<T_{\min X}$ or $\sigma_{\max
X}(K)<T_{\max X}$. Suppose without loss of generality that
$\sigma_{\min X}(K)<T_{\min X}$. Then $\sigma_{m}(K^g)<T_{m}$
where $m\in\{1,\ldots,k\}$ is such that $T_{\min X}^g=T_m$. As $X^g$
covers $T_m$, it follows that $\{K,X\}^g=\{K^g,X^g\}\in
E_1$. Now let $\{K,X\}\in E_2$. Then $X$ is involved in $K$ and hence 
$X^g$ is involved in $K^g$, whence $\{K^g,X^g\}\in E_2$. Thus $G_\omega$
preserves adjacency in $\Gamma(G,\E)$. 
Moreover, under the conjugation action of $G_\omega$, the sets $\K_\omega(\E)$
and $\X$ are $G_\omega$-orbits. We claim that $E_1$ and $E_2$ are also
$G_\omega$-orbits. Suppose that 
$\{K_1,X_1\},\ \{K_2,X_2\}\in E_1$. There there exist $i_1,\ i_2$
such that $T_{i_1}\in\supp X_1$, $T_{i_2}\in\supp X_2$,
$\sigma_{i_1}(K_1)<T_{i_1}$, and $\sigma_{i_2}(K_2)<T_{i_2}$. Since
$G_\omega$ is transitive on $T_1,\ldots,T_k$, there is
an element $g\in G_\omega$ such that $T_{i_1}^g=T_{i_2}$. Then
$T_{i_2}\in\supp X_1^g\cap\supp X_2$, and so
Theorem~\ref{stripth} implies that $X_1^g=X_2$. We also have that
$\sigma_{i_1}(K_1)^g=\sigma_{i_2}(K_1^g)<T_{i_2}$. Since $K_2$ is the
unique element in $\K_\omega(\E)$ with proper projection in $T_{i_2}$
we have $K_1^g=K_2$. Thus $\{K_1,X_1\}^g=\{K_2,X_2\}$, and so
$G_\omega$ is transitive on $E_1$. 
Now let $\{K_1,X_1\},\ \{K_2,X_2\}\in
E_2$. Then $X_1$ is involved in $K_1$ and $X_2$ is involved in $K_2$.
There is an element $g\in G_\omega$ such that $X_1^g=X_2$, which
implies that $K_1^g=K_2$. Thus $\{K_1,X_1\}^g=\{K_2,X_2\}$, and
$G_\omega$ is transitive on $E_2$. 
Finally suppose that the elements of $\X$ have
$E_1$-valency~2 and let $\{K,X\}\in E_1$. Suppose without loss of
generality that $\sigma_{\min X}(K)<T_{\min X}$ and let $g\in
(G_\omega)_{\{K,X\}}$. Then $\{T_{\min X},T_{\max X}\}^g=\{T_{\min
X},T_{\max X}\}$, and so either $$
T_{\min X}^g=T_{\min X}\quad\mbox{and}\quad T_{\max X}^g=T_{\max X}
$$
or 
$$
T_{\min X}^g=T_{\max X}\quad\mbox{and}\quad T_{\max X}^g=T_{\min X}
$$
In the latter case we would have $\sigma_{\min X}(K)^g=\sigma_{\max
X}(K)$. As $X$ has $E_1$-valency 2, we have that $\sigma_{\max
X}(L)<T_{\max X}$ for a unique 
$L\in\K_\omega(\E)$ and this $L$ is different from $K$, which is a
contradiction. Thus $T_{\min X}^g=T_{\min X}$ and $T_{\max
X}^g=T_{\max X}$ must hold.
Therefore $g\in \norm{G_\omega}{T_{\min
X}}\cap\norm{G_\omega}{T_{\max X}}$.
Conversely suppose that $g\in \norm{G_\omega}{T_{\min
X}}\cap\norm{G_\omega}{T_{\max X}}$. Then clearly $g\in\norm
{G_\omega}X$. Moreover, $\sigma_{\min X}(K)^g=\sigma_{\min
X}(K^g)$, and since $K$ is the unique element of $\K_\omega(\E)$ such that
$\sigma_{\min X}(K)<T_{\min X}$, it follows that $K^g=K$. Therefore
$\{K,X\}^g=\{K,X\}$. Hence property~\textsf{Prop1S[c]} holds.

\textbf{\textsf{Prop1S[d]}}\quad
By Proposition~\ref{6.1th} and Lemma~\ref{uniqueext} the groups $T$, $A_i$ 
(which is $G_\omega$-conjugate to $B_i$), and $\sigma_{\supp X_i}(M_{
\overline\omega})$ are as in Table~\ref{isomfact}. 
By Lemma~\ref{uniqueext}, the group $\overline M_\omega$ is a direct product of
pairwise disjoint strips, and each such strip $Y$ is self-normalising
in $\sigma_{\supp Y}(M)$ (see Lemma~\ref{normstriplemma}). Thus
$\overline M_\omega$ is a self-normalising subgroup of
$M$. Hence~\cite[Theorem~4.2A]{dm} implies that
$\cent{\sym\overline\Omega}M=1$. Thus $\overline G$ can be embedded
into $\aut M$, and so $\overline G$ is  quasiprimitive and $M$ is its unique
minimal normal subgroup.

By Lemma~\ref{knorm}, if one of the rows~1--3 of Table~\ref{isomfact} is valid, 
then $K_i=\overline K_i$ for all $i$, and so $M_\omega=\overline M_\omega$. 
Thus the sets $\overline\Omega$ and $\Omega$ can be identified naturally, and
the groups $G$ and $\overline G$ are naturally permutationally isomorphic.

If $T$ is as in row~4 of Table~\ref{isomfact}, then it follows from 
Lemma~\ref{uniqueext}
that $\norm M{M_\omega}/M_\omega$ is an elementary 
abelian $2$-group with rank at most $k/2$. Thus each block in $\overline\Omega$
has size dividing $k/2$. Further, it follows from~\cite[Theorem~4.2A]{dm}
that $\cent{\sym\Omega}M$ is also an elementary abelian 2-group of rank at 
most $k/2$. As $M$ is a minimal normal subgroup of $G$ and is faithful on 
$\Omega$, we obtain that $N\cap M=1$, and so $N\leq\cent{\sym\Omega}M$. 
Thus $N$ an elementary abelian 2-group of rank at most $k/2$.\hfill$\Box$

\subsection{A converse of Theorem~\ref{main1s}}\label{c1s}\quad
Theorem~\ref{main1s} can be reversed in the following sense. Suppose that $G$
is an innately transitive  group on $\Omega$ with non-abelian plinth $M$,
and let $T_1,\ldots,T_k$ be the simple normal subgroups of
$M$. Let $\omega\in\Omega$. Assume that
$M_\omega$ is a direct product of pairwise disjoint strips, each of length
$2$.
Let $\Y$ denote the
set of strips involved in $M_\omega$, say $\Y=\{Y_1,\ldots,Y_{k/2}\}$.
Let $X_1$, $A_1$, $B_1$ be subgroups
of $M$, 
and let $\Gamma=(V\cup\mathcal Y,E_1\cup E_2)$ be a bipartite graph 
satisfying properties {\sf Prop1S[b]} and {\sf Prop1S[c]}, that is, the 
following all hold:
\begin{enumerate}
\item[(i)] $X_1$ is a full strip of length $2$;
\item[(ii)] $A_1$ is a proper subgroup of $T_{\min X_1}$ and 
$B_1$ is a proper subgroup of $T_{\max X_1}$;
\item[(iii)] $A_1$ and $B_1$ are
conjugate under $G_\omega$;
\item[(iv)] $X_1(A_1\times B_1)=T_{\min X_1}\times T_{\max X_1}$,   
$X_1\cap(A_1\times B_1)=\sigma_{\supp X_1}(M_{\omega})$;
\item[(v)] $\norm {G_\omega}{T_{\min X_1}\times T_{\max X_1}}=\norm{
G_\omega}{X_1}=\norm{G_\omega}{A\times B}$;
\end{enumerate}
and also the group $G_\omega$ induces a group of automorphisms of
the bipartite graph $\Gamma$, such that $V$, $\Y$, $E_1$, 
and $E_2$ are $G_\omega$-orbits. 
Further, each element of $\Y$ is adjacent to one edge 
or two edges from $E_1$, and
 one edge of $E_2$.
If the elements of $\Y$ are adjacent to two elements of
$E_1$ then the following must also hold: if $\{v,Y\}\in E_1$, 
for some $v\in V$ and $Y\in\Y$, then $(G_\omega)_{\{v,Y\}}=\norm{G_\omega}
{T_{\min Y}}\cap \norm{G_\omega}{T_{\max Y}}$.

We claim that $\supp X_1=\supp Y_j$ for some $j\in\{1,\ldots,k/2\}$. If this is not true
then 
we have  $\sigma_{\supp X_1}(M_\omega)=\sigma_{T_{\min X_1}}(M_\omega)\times \sigma_{T_{\max X_1}}(M_\omega)$, and by property
(iv) above, this is contained in $X_1\cap(A_1\times B_1)$, which is a
contradiction. Hence we may assume, without loss 
of generality, that $\supp X_1 = \supp Y_1$. Let $v_1,\ v_2\in V$ be
such that $\{v_1,Y_1\}\in E_1$ and $\{v_2,Y_1\}\in E_2$, and  define
$$
K_{v_1,Y_1}=\left\{\begin{array}{ll}
A_1\times T_{\max Y_1} & \mbox{if $Y_1$ has $E_1$-valency~$2$;}\\
A_1\times B_1 & \mbox{otherwise.}
\end{array}\right.
$$  
For $v\in V$ and $Y\in\Y$, define
$$
K_{v,Y}=
\left\{\begin{array}{ll}
(K_{v_1,Y_1})^g & \mbox{if
$\{v,Y\}\in E_1$, and $g\in G_\omega$ is such that $\{v_1,Y_1\}^g=\{v,Y\}$;}\\
X_1^g & \mbox{if $\{v,Y\}\in E_2$, and $g\in G_\omega$ is such that $\{v_2,Y_1\}^g=\{v,Y\}$;}\\
T_{\min Y}\times T_{\max Y} & \mbox{otherwise.}\end{array}\right.
$$
We claim that the definition of $K_{v,Y}$ is independent of the chosen
$g\in G_\omega$. Indeed, if two elements  $g_1,\ g_2\in G_\omega$ are such
that $\{v_1,Y_1\}^{g_1}=\{v_1,Y_1\}^{g_2}=\{v,Y\}$, then
$g_1g_2^{-1}\in(G_\omega)_{\{v_1,Y_1\}}$. If  $Y_1$ has $E_1$-valency~2
then by
assumption, 
$$
g_1g_2^{-1}\in\norm{G_\omega}{T_{\min
Y_1}}\cap\norm{G_\omega}{T_{\max
Y_1}}=\norm{G_\omega}{A_1}\cap\norm{G_\omega}{B_1}.
$$
Hence
$$
(K_{v_1,Y_1})^{g_1g_2^{-1}}=(A_1\times T_{\max
Y_1})^{g_1g_2^{-1}}=A_1^{g_1g_2^{-1}}\times (T_{\max
Y_1})^{g_1g_2^{-1}}=A_1\times T_{\max Y_1}=K_{v_1,Y_1}
$$
If $Y_1$ has $E_1$-valency~1 then, as $g_1g_2^{-1}\in\norm{G_\omega}
{T_{\min Y_1}\times T_{\max Y_1}}$, property~(v) implies that
$$
(K_{v_1,Y_1})^{g_1g_2^{-1}}=(A_1\times B_1)^{g_1g_2^{-1}}=A_1\times
B_1=K_{v_1,Y_1}.
$$
Thus $(K_{v_1,Y_1})^{g_1}=(K_{v_1,Y_1})^{g_2}$, as claimed. 
Similarly if $\{v_2,Y_1\}^{g_1}=\{v_2,Y_1\}^{g_2}$ 
for some $g_1,\ g_2\in G_\omega$
then 
$$
g_1g_2^{-1}\in\norm{G_\omega}{T_{\min Y_1}\times T_{\max
Y_1}}=\norm{G_\omega}{X_1}
$$
and so $X_1^{g_1}=X_1^{g_2}$. So also in this case the definition of
$K_{v,Y}$ is independent of the chosen element $g$.

We now claim
that for
each $v\in V$, $Y\in\Y$, and $g\in G_\omega$ we have 
\begin{equation}\label{rule}
K_{v^g,Y^g}=(K_{v,Y})^g.
\end{equation}
Indeed suppose that $\{v,Y\}\in E_1$. Then there is an element $g_1\in 
G_\omega$ such that $\{v_1,Y_1\}^{g_1}=\{v,Y\}$, and so
$\{v^g,Y^g\}=\{v_1,Y_1\}^{g_1g}$. 
Now $K_{v^g,Y^g}$ was defined above as $(K_{v_1,Y_1})^{g_1g}$, and $K_{v,Y}$
was defined as $(K_{v_1,Y_1})^{g_1}$. Thus (\ref{rule}) holds in this case.
Similarly, if  $\{v,Y\}\in E_2$ then there is an element $g_1\in G_\omega$ 
such that $\{v_2,Y_1\}^g=\{v,Y\}$, and the same argument shows 
that~\eqref{rule} holds. Finally, if
$\{v,Y\}\not\in E_1\cup E_2$ then $\{v^g,Y^g\}\not\in E_1\cup
E_2$. Hence $K_{v,Y}=T_{\min Y}\times T_{\max Y}$ and
$K_{v^g,Y^g}=T_{\min Y^g}\times T_{\max Y^g}$. As $\supp Y^g=(\supp
Y)^g$ we have that $(T_{\min Y}\times T_{\max Y})^g=T_{\min
Y^g}\times T_{\max Y^g}$. Hence \eqref{rule} holds in all cases.

Now we note that $K_{v,Y}$ is a subgroup of $T_{\min Y}\times
T_{\max Y}$, and, as the elements of $\Y$ are pairwise disjoint
strips, we can define 
$$
K_v=\prod_{Y\in \Y}K_{v,Y}
$$
and set $\K=\{K_v\ |\ v\in V\}$. 

Our next task is to prove that $\K$ is a Cartesian system, that is,
equations~\eqref{csdef1} and~\eqref{csdef2} hold. To help ourselves
with this, first we prove analogous properties for the subgroups
$K_{v,Y}$ of the $K_v$. 
If $Y_1$ has $E_1$-valency 1 then, $Y_1$ is adjacent to two vertices
$v_1$ and $v_2$, and, by property (iv), we have
\begin{equation}\label{prod1}
K_{v_1,Y_1}K_{v_2,Y_1}=(A_1\times B_1)X_1=T_{\min Y_1}\times T_{\max
Y_1}
\end{equation}
and
\begin{equation}\label{int1}
K_{v_1,Y_1}\cap K_{v_2,Y_1}=(A_1\times B_1)\cap X_1=\sigma_{\supp
Y_1}(M_\omega).
\end{equation}
Suppose now that $Y_1$ has $E_1$-valency 2, and let $v_3$ be a vertex
such that $v_1\neq v_3$ and $\{v_3,Y_1\}\in E_1$. 
Then there is some element $g\in
G_\omega$ such that $A^g=B$. Then we must have that $T_{\min
Y_1}^g=T_{\max Y_1}$, and so, by the conditions above, $g$ must
interchange $v_1$ and $v_3$. Then the argument above shows that $K_{v_3,Y_1}=T_{\min
Y_1}\times B$. Thus, by (iv),
\begin{equation}\label{prod2}
\{K_{v_1,Y_1},K_{v_2,Y_1},K_{v_3,Y_1}\}\mbox{ is a
strong multiple factorisation of $T_{\min Y_1}\times T_{\max Y_1}$}
\end{equation}
and 
\begin{equation}\label{int2}
K_{v_1,Y_1}\cap K_{v_2,Y_1}\cap K_{v_3,Y_1}= \sigma_{\supp
Y_1}(M_\omega).
\end{equation}

Now we are ready to show that~\eqref{csdef1} and~\eqref{csdef2} hold. 
As $M_\omega$ and
the elements of $\K$ are direct products of their projections under
$\sigma_{\supp Y}$, for $Y\in\Y$, it suffices to prove that 
$$
\bigcap_{v\in V}K_{v,Y}=\sigma_{\supp Y}(M_\omega)\quad\mbox{and}\quad
K_{v,Y}\left(\bigcap_{v'\neq v}K_{v',Y}\right)=T_{\min Y}\times
T_{\max Y}
$$
holds for all $Y\in\Y$. If $g\in G_\omega$ such that $Y_1^g=Y$ then,
by~\eqref{int1} and~\eqref{int2}, we
have that
$$
\bigcap_{v\in V}K_{v,Y}=X_1^g\cap(A_1\times B_1)^g=\left(X_1\cap (A_1\times 
B_1\right))^g=\sigma_{\supp Y_1}(M_\omega)^g=\sigma_{\supp
Y}(M_\omega).
$$
Also, using, \eqref{prod1} and~\eqref{prod2},
$$
K_{v,Y}\left(\bigcap_{v'\neq
v}K_{v',Y}\right)=\left(K_{v^{g^{-1}},Y_1}\left(\bigcap_{v'\neq
v}K_{v'^{g^{-1}},Y_1}\right)\right)^g=\left(T_{\min Y_1}\times T_{\max
Y_1}\right)^g=T_{\min Y}\times T_{\max Y}.
$$
Hence~\eqref{csdef1} and~\eqref{csdef2} hold. 

It remains to prove that $\K$ is $G_\omega$-invariant. Let $v\in V$
and $g\in G_\omega$. Then
$$
K_v^g=\left(\prod_{Y\in\Y}K_{v,Y}\right)^g=\prod_{Y\in\Y}K_{v^g,Y^g}=K_{v^g}.
$$
Thus $K_v^g\in \K$. Therefore $\K$ is a $G_\omega$-invariant Cartesian
system of subgroups in $M$. It also follows from the last displayed
equation that the actions of $G_\omega$  on $\K$ and on $V$ are equivalent, 
and so
$G_\omega$ is transitive on $\K$. It follows from the
definition of the $K_v$ that $\E(\K)\in\cd{1S}G$ and that $\Gamma=\Gamma(G,\E)$. 

\medskip

Suppose that $G$ is an innately transitive group with a non-abelian plinth 
$M$ and a point stabiliser $M_\omega$ is a direct product of pairwise 
disjoint strips with length~2 such that
the isomorphism types of these groups are as prescribed by Table~\ref{isomfact}.
Then, as in the previous sections, it is sometimes possible to describe the
elements of $\cd{1S}G$ via studying the action of $G_\omega$ on the set of
strips in $M_\omega$. This phenomenon is illustrated in the following example.

\begin{example}
Suppose that $\Gamma=(V_1\cup V_2,E_1\cup E_2)$ is a bipartite graph 
with an automorphism group $A$ satisfying
the Combinatorial Property where $V_1$, $V_2$, and $A$ play the r\^ole of 
$\K_\omega(\E)$, $\X$, and $G_\omega$, respectively. 
Suppose that $V_2$ has 4~vertices and 
$A$ induces a group isomorphic to 
$\alt 4$ on $V_2$. As each vertex in $V_2$ is adjacent to exactly 
one edge in $E_2$ it follows that $V_1$ must have~2 or~4 vertices. If $V_1$ has
2 vertices, $v_1$ and $v_2$ say, then the set of vertices in $V_2$ 
connected to 
$v_1$ via $E_2$ is a block for the action of $A$. Such a block would have~2 elements, and this is impossible, as $\alt 4$ has no non-trivial blocks. 
Hence $|V_1|=4$. 
If a vertex of $V_2$ is adjacent with
2~edges in $E_1$ then $E_1$ must have 8~elements. Thus $E_1$ cannot be an
$A$-orbit, as $|A|$ is not divisible by 8. Thus each element of $V_2$ must be adjacent with exactly one edge of $E_1$. Suppose without loss of generality 
that $\{v_1,u_1\}\in E_1$ and $\{v_2,u_1\}\in E_2$ for some $u_1\in V_2$. Then 
$v_1\neq v_2$, and so $A_{u_1}\leq A_{v_1}\cap A_{v_2}=1$. This is a 
contradiction since $|A:A_{u_1}|=|V_1|=4$. Thus $E_2$ cannot have $4$ 
elements,  and so no such graph $\Gamma$ exists.

This simple graph theoretic argument shows that if $G$ is an innately transitive group
with plinth $M=T_1\times\cdots\times T_8$ such that a point
stabiliser $M_\omega$ is the direct product of 4~pairwise disjoint strips
and $G_\omega$ induces a group permutationally 
isomorphic to $\alt 4$ on these strips then $\cd{1S}G=\emptyset$.
\end{example}


\begin{thebibliography}{BamP04}

\bibitem[BP98]{bad:fact}
Robert~W. Baddeley and Cheryl~E. Praeger,
\newblock On classifying all full factorisations and multiple-factorisations of
  the finite almost simple groups,
\newblock {\em J. Algebra}, 204(1):129--187, 1998.


\bibitem[BP03]{bad:quasi}
R.~W. Baddeley and C.~E. Praeger,
\newblock On primitive overgroups of quasiprimitive permutation groups,
\newblock {\em J. Algebra}, 263(2):294--344, 2003.



\bibitem[BPS04]{recog}
Robert W. Baddeley, Cheryl E. Praeger and Csaba Schneider. 
\newblock Transitive simple subgroups of wreath products in product
action. {\em J. Austral. Math. Soc.} 77(1):55-72, 2004.

\bibitem[BPSxx]{transcs}
\newblock Robert~W. Baddeley, Cheryl~E. Praeger, and Csaba Schneider.
Innately transitive subgroups of wreath products in product action.  
To appear in {\em Trans. Amer. Math. Soc.} arXiv.org/math.GR/0312352.

\bibitem[PSxx]{intrans} Robert W. Baddeley, Cheryl E. Praeger, 
and Csaba Schneider. 
\newblock Intransitive Cartesian decompositions preserved by innately
transitive groups. Submitted.
\newblock arxiv.org/math.GR/0405241.


\bibitem[BamP04]{bp}
\newblock John Bamberg and Cheryl E. Praeger. Finite permutation groups with
a transitive minimal normal subgroup. {\em
  Proc. London. Math. Soc. (3)} 89(1):71-103, 2004.


\bibitem[Atlas]{atlas}
J.~H. Conway, R.~T. Curtis, S.~P. Norton, R.~A. Parker, and R.~A. Wilson.
\newblock {\em Atlas of finite groups}.
\newblock Oxford University Press, Oxford, 1985.


\bibitem[DM96]{dm}
John~D. Dixon and Brian Mortimer,
\newblock {\em Permutation groups},
\newblock Springer-Verlag, New York, 1996.

\bibitem[Kle87]{kleidman}
Peter~B. Kleidman.
\newblock The maximal subgroups of the finite $8$-dimensional orthogonal groups
  ${P}\Omega^+_8(q)$ and of their automorphism groups.
\newblock {\em J. Algebra}, 110(1):173--242, 1987.

\bibitem[Kov89a]{kov:blowups}
L.~G. Kov{\'a}cs.
\newblock Primitive subgroups of wreath products in product action.
\newblock {\em Proc. London Math. Soc. (3)}, 58(2):306--322, 1989.

\bibitem[Kov89b]{kov:decomp}
L.~G. Kov{\'a}cs.
\newblock Wreath decompositions of finite permutation groups.
\newblock {\em Bull. Austral. Math. Soc.}, 40(2):255--279, 1989.


\bibitem[LPS90]{lps:max}
Martin~W. Liebeck, Cheryl~E. Praeger, and Jan Saxl.
\newblock The maximal factorizations of the finite simple groups and their
  automorphism groups.
\newblock {\em Mem. Amer. Math. Soc.}, 86(432):iv+151, 1990.



\bibitem[Pra93]{prae:quasi}
Cheryl~E. Praeger.
\newblock An {O}'{N}an-{S}cott theorem for finite quasiprimitive permutation
  groups and an application to $2$-arc transitive graphs.
\newblock {\em J. London Math. Soc. (2)}, 47(2):227--239, 1993.


\bibitem[PS02]{charfact} Cheryl E. Praeger and Csaba Schneider. 
\newblock Factorisations of characteristically simple groups.
\newblock {\em J. Algebra}, 255(1):198--220, 2002.




\bibitem[PS03]{design}
Cheryl~E. Praeger and Csaba Schneider,
\newblock Ordered triple designs and wreath products of groups.
\newblock In Darlene R. Goldstein (Ed.) {\em Science and Statistics: A
Festschrift for Terry Speed}. Institute of Mathematical Statistics,
Lecture Notes -- Monograph Series, volume {\bf 40}, pages 103--113,
2003.

\end{thebibliography}
\end{document}